\documentclass{article}

\usepackage[all]{nowidow}
\usepackage[protrusion=true, final]{microtype}
\usepackage{josealvim}
\usepackage{hyperref}
\usepackage[backend=biber]{biblatex}
\usepackage{quiver}
\usepackage[inline]{enumitem}
\usepackage{import}
\usepackage{verbatim}
\usepackage[utf8]{inputenc}
\begin{document}
	\title{
		-- \(\QSet\) \textit{\&} Friends --
		\\
        \Large
        Categorical Constructions and Categorical Properties
	}
	\author{
        José G. Alvim\email{alvim@ime.usp.br}
		\and 
        Hugo L. Mariano\email{hugomar@ime.usp.br}
		\and
        Caio de A. Mendes\email{caio.mendes@usp.br}
	}
	\begin{titlepage}
	
		\maketitle
		\begin{abstract}
            This work mainly concerns the -here introduced- category of \(\quantale Q\)-sets 
            and functional morphisms, where \(\quantale Q\) is a commutative semicartesian quantale.
            We describe, in detail, the limits and colimits of this complete and
            cocomplete category and prove that it has a classifier for regular 
            subobjects. Moreover, we prove that it is \(\kappa^+\)-locally 
            presentable category, where \(\kappa=max\{|\quantale Q|, \aleph_0)\}\) 
            and describe a hierarchy of semicartesian monoidal closed structures 
            in this category. Finally, we discuss the issue of “change of basis” 
            induced by appropriate morphisms between the parametrizing quantales
            involved in the definition of \(\quantale Q\)-sets. 
            In a future work we will address such questions in the full 
            subcategory given by all Scott-complete \(\quantale Q\)-sets (see \cite{AMM2023completeness}).
		\end{abstract}       
		\thispagestyle{empty}
	\end{titlepage}

	\tableofcontents
    \newpage
    
	\section{Introduction}


\subsection*{History \textit{\&} Motivation}









Sheaf Theory is a well established area of research with applications
in Algebraic Topology \cite{dimca2004sheaves}, Algebraic Geometry 
\cite{grothendieck1972topos}, Geometry \cite{kashiwara2013sheaves}, 
Logic \cite{maclane1992sheaves}, and others.  

A sheaf on a topological space \(X\) is a functor 
\(F:(\mathcal{O}(X),\subseteq)^{op}\to\Set\) that satisfies certain 
gluing properties expressed by an equalizer diagram, where  
\((\mathcal{O}(X),\subseteq)\) denotes the category whose objects are 
the elements of the set  \(\mathcal{O}(X)\) of all open subspaces of 
\(X\), and the morphisms are given by set inclusions. There are other 
equivalent ways to express this definition, but the diagrammatic 
approach makes clear that the elements/points of \(X\) are not 
necessary. Thus it is straightforward to define sheaves for “spaces 
without points”, that is, for a category \(H,\leq)\) for a given 
locale \(H\). The category of sheaves on locales was studied, for 
example, by Borceux in \cite{borceux1994handbook3}. 

In the 1970s, the topos of sheaves over a locale/complete Heyting 
algebra \(\mathbb{H}\), denoted as \(\Sh(\mathbb{H})\), was described,
alternatively, as a category of \(\mathbb{H}\)-sets 
\cite{fourman1979sheaves}. More precisely, in 
\cite{borceux1994handbook3}, there were three categories whose objects 
were locale valued sets that are equivalent to the category of sheaves 
over a locale \(\mathbb{H}\). 



Nevertheless, there is a non-commutative and non-idempotent 
generalization of locales called “quantales”, introduced by C.J. Mulvey 
\cite{mulvey86quantales}. Quantales show up in logic 
\cite{yetter_1990}, and in the study of \(C^*\)-algebras 
\cite{rosenthal1990quantales}. 

Later, more general categories have been proposed, replacing locales by 
the Mulvey's quantales 
studied in \cite{coniglio2000non}. Instead considering the traditional 
idempotent non-commutative quantales, that arose from certain 
\(C^*\)-algebras and its relationships with quantum physics; we are 
following proposals like in \cite{HOHLE199115} and \cite{zambrano}, 
that have connections with affine, fuzzy, and continuous logic. In this 
work we consider a class of commutative and integral/semicartesian 
quantales, which includes both the quantales of the ideals of 
commutative unital rings, MV-Algebras, Heyting Algebras and 
\(([0,1], \leq, \cdot)\) -- which is isomorphic to 
\(({[0,\infty]}, \geq, +)\). 

So far as we know, there are three notions of sheaves on 
\textit{right-sided and idempotent} quantales: in 
\cite{borceux1986quantales}, sheaves on  quantales are defined with the 
goal of forming Grothendieck toposes from quantales. 
In \cite{miraglia1998sheaves}, the sheaf definition  preserves an 
intimate relation with \(\quantale Q\)-sets, an object introduced in 
the paper as a proposal to generalize \(\Omega\)-sets, defined in 
\cite{fourman1979sheaves}, for \(\Omega\) a complete Heyting algebra%
\footnote{
    Given a proper notion of morphisms between \(\Omega\)-sets, the 
    resulting category is equivalent to the category of sheaves on 
    \(\Omega\).
}. 

More recently, in \cite{aguilar2008sheaves} and \cite{tenorio2022sheaves}, sheaves are functors that 
make a certain diagram an equalizer. Herein we study sheaves on 
\emph{semicartesian} quantales. Our approach is similar to the last one 
but, since every idempotent semicartesian quantale is a locale 
(Proposition \ref{prop:idempotent + semicartesian = locale}), our 
axioms and theirs are orthogonal in some sense.

Besides, there is an extensive work about sheaves on 
\textit{involutive quantale}, which goes back to ideas of Bob Walters 
\cite{walters1981sheaves} -- which were recently studied by Hans 
Heymans, Isar Stubbe \cite{heymans2012grothendieck}, and Pedro Resende 
\cite{resende2012groupoid} -- for instance.



               This work mainly concerns the -here introduced- category of \(\quantale Q\)-sets 
            and functional morphisms, where \(\quantale Q\) is a commutative semicartesian quantale.
             In a future work we will address such questions in the full 
            subcategory given by all Scott-complete \(\quantale Q\)-sets (see \cite{AMM2023completeness}) .

\subsection*{Main results and the paper's structure}

\begin{enumerate}
    \item We describe, in detail, the limits and colimits of this complete and cocomplete category; 
    \item  We describe generators and prove that it has a classifier for regular subobjects;
         \item We prove that it is \(\kappa^+\)-locally  presentable category, where \(\kappa=max\{|\quantale Q|, \aleph_0)\}\); 
    \item We describe a hierarchy of semicartesian monoidal closed structures   in this category;
    \item  We discuss the issue of “change of basis”  induced by appropriate morphisms between the parametrizing quantales involved in the definition of \(\quantale Q\)-sets.

\end{enumerate}





    \section{Preliminaries}

\subsection{Quantales}
\label{sec:quantales}

\begin{definition}
    A \emph{quantale} is a type of structure 
    \(\quantale Q = (|\quantale Q|, \leq, \tensor)\) for which 
    \((|\quantale Q|, \leq)\) is a complete lattice; 
    \((|\quantale Q|, \tensor)\) is a semigroup%
    \footnote{
        \lat{i.e.} the  binary operation 
        \(\tensor:\quantale Q\times\quantale Q\to\quantale Q\) 
        (called multiplication) is associative.
    }; and, moreover, \(\quantale Q\) is required to satisfy the 
    following distributive laws: for all \(a\in\quantale Q\) and 
    \(B\subseteq\quantale Q\),
    \begin{align*}
            a \tensor \left(\bigvee_{b\in B} b\right) 
        &=	\bigvee_{b\in B}\left(a\tensor b\right)
        \\
            \left(\bigvee_{b\in B} b\right) \tensor a 
        &=	\bigvee_{b\in B}\left(b\tensor a\right)
    \end{align*}
    We denote by \(\extent\quantale Q\) the subset of 
    \(\quantale Q\) comprised of its idempotent elements.
\end{definition}

\begin{remark}~\label{remark:quantale properties}
    \begin{enumerate}
        \item 
            In any quantale \(\quantale Q\) the 
            multiplication is increasing in both entries;
        \item 
            Since \(\bot\) is also the supremum of 
            \(\emptyset\), for any \(a\),
            \(a\tensor\bot=\bot=\bot\tensor a\)
        \item 
            Since \(\top\) is \(\sup\quantale Q\), then
            \(
                \top\tensor\top=
                \sup_{a,b}a\tensor b =
                \sup\img\tensor
            \)
    \end{enumerate}
\end{remark}

\begin{remark} 
    If \((\quantale Q, \leq)\) is a complete lattice for which the 
    binary infimum satisfies the above distributive laws, the 
    resulting quantale has \(\top\) as its unit and is -- in fact 
    -- a locale. Conversely, every locale is a unital quantale in
    such a manner.
\end{remark}

\begin{definition}
    A quantale \(\quantale Q\) is said to be
    \begin{enumerate}[label=\(\bullet\)]    
        \item \emph{bidivisible}
            when \[
                a\leq b\implies 
                \exists \lambda,\rho: a\tensor\rho = b = \lambda\tensor a
            \]
            left (right) divisibility means to drop the 
            \(\rho\) (\(\lambda\)) portion of the axiom.
        \item \emph{integral}      
            when \(\top\tensor a = a = a\tensor\top\). 
            we say it's right-sided when the right equality
            holds, and left-sided when the left equality 
            holds.
        \item \emph{unital}
            when \(\tensor\) has a unit;
        \item \emph{semicartesian} 
            when \(a \tensor b \leq a\wedge b\)
        \item \emph{commutative}
            when \(\tensor\) is;
        \item \emph{idempotent}
            when \(a \tensor a = a\);
        \item \emph{linear and strict}
            when \(\leq\) is a linear order and for \[
                a\neq \bot\ and \ a\tensor b = a\tensor c\ or\ b\tensor a = c\tensor a  \implies 
                b = c
            \]
        \item \emph{strong} 
            when for any \(e\) and \(A\)
            \cite[cf. p. 30]{HOHLE199115},
        \[
            e=e\tensor e \implies 
            e\leq\bigvee_{a\in A} a \implies 
            e\leq\bigvee_{a\in A} a \tensor a
        \]
    \end{enumerate}
    We offer the following diagram to explain some of 
    the relations between those definitions:
    
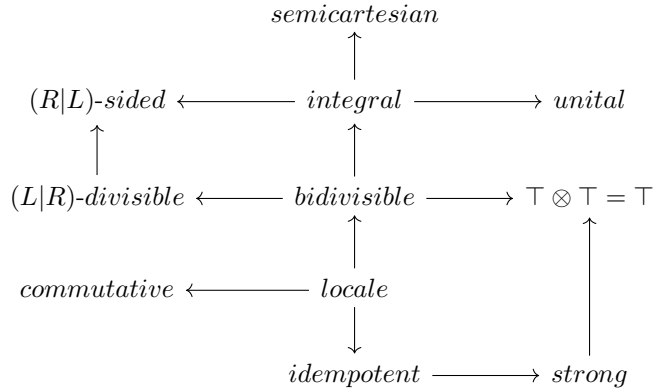
\[\begin{tikzcd}
        & semicartesian \\
        {(R|L)\text-sided} & integral & unital \\
        {(L|R)\text-divisible} & bidivisible & {\top\tensor\top=\top} \\
        commutative & locale \\
        & idempotent & strong
        \arrow[from=3-2, to=3-1]
        \arrow[from=3-2, to=2-2]
        \arrow[from=4-2, to=3-2]
        \arrow[from=4-2, to=4-1]
        \arrow[from=2-2, to=2-1]
        \arrow[from=2-2, to=1-2]
        \arrow[from=2-2, to=2-3]
        \arrow[from=4-2, to=5-2]
        \arrow[from=5-2, to=5-3]
        \arrow[from=3-1, to=2-1]
        \arrow[from=5-3, to=3-3]
        \arrow[from=3-2, to=3-3]
    \end{tikzcd}\]
\end{definition}

\begin{example}
    Locales are -- perhaps -- the best example of quantales that 
    are commutative, idempotent, integral (and hence both 
    semicartesian and right-sided), divisible and strong (both 
    trivially). Among which, and of special significance to 
    Sheaf Theory, are the locales of open subsets of a topological 
    space \(X\), where the order relation is given by the
    inclusion, the supremum is the union, and the finitary infimum 
    is the intersection.
\end{example}

\begin{example}
    The extended half-line \([0,\infty]\) with order the inverse 
    order -- \(\geq\) --, and the usual sum of real numbers as the 
    monoidal operation. 

    Since the order relation is \(\geq\), the top element is \(0\)
    and the bottom elements is \(\infty\). We call this the Lawvere
    quantale due to its relation to Lawvere spaces (related to 
    metric spaces).

    An isormophic quantale is the quantale given by \([0,1]\) this 
    time with the usual order -- and multiplication as the monoidal
    operation. The isomorphism is given by 
    \[
        (x\in[0,1]) \mapsto (-\ln(x)\in[0,\infty])
    \]

    A sub-example is the extended natural numbers 
    \(\mathbb{N}\cup\{\infty\}\), as a restriction of the the 
    Lawvere quantale, which is related to distance on graphs.

    All of these are unital but not locales. In fact, they are 
    integral too. 
\end{example}

\begin{remark}
    When we claim that a quantale isn't a locale, we do not mean to
    say that the underlying poset isn't a locale. Merely that the 
    monoidal product isn't the poset's meet operation. 
\end{remark}

\begin{example}\label{example: unital quantales}
    We list below some more examples of unital quantales that are 
    not locales:
    \begin{enumerate}
        \item 
            The set \(\mathcal{I}(R)\) of ideals of a 
            commutative and unital ring \(R\) with order 
            \(\subseteq\), and the multiplication as the 
            multiplication of ideals. The supremum is the 
            sum of ideals, the top element is \(R\) and the 
            trivial ideal is the bottom;
        \item 
            The set \(\mathcal{RI}(R)\) of right (or left) 
            ideals of an unital ring \(R\) with the same 
            order and multiplication of the above example. 
            Then the supremum and the top and the bottom 
            elements are also the same of 
            \(\mathcal{I}(R)\);
        \item 
            The set of closed right (or left) ideals of a 
            unital \(C^*\)-algebra, the order is the 
            inclusion of closed right (or left) ideals, and 
            the multiplication is the topological closure 
            of the multiplication of the ideals.
    \end{enumerate}
    For more details and examples we recommend 
    \cite{rosenthal1990quantales}. 
\end{example}

\begin{example}[The “slice” construction]
    Given a quantale \(\quantale Q\), we can form 
    intervals between two extremes $a \leq b$:
    \[
        [a, b] = \set{p\in\quantale Q}{a\leq p\leq b}
    \]
    This is evidently a complete lattice and the inclusion \([a, b] \hookrightarrow {\quantale Q}\) preserves non empty sups and non empty infs.
    
    Now, suppose \(\quantale Q\) is (commutative) semicartesian and unital, and consider 
    -- given any \(e\in Idem(\quantale Q)\) and \(b\in \quantale Q\) such that \(e \leq b\) -- the set \([e,b]\) is such that

    \(e \leq x,y\leq b \implies e = e \tensor e \leq x\tensor y  \leq b \tensor b \leq b \otimes \top = b\) 
    
    hence that \(\tensor\restriction[e,b]\) is a subsemigroup. Moreover, that this  has the structure of 
    a quantale \(x \tensor \bigvee_i y_i = \bigvee_i x \tensor y_i\): there is not to check concerning the distributivity of \(x \in [e,b]\)  with non-empty sups; for empty sups, we just have to notice that 
    
    \( e = e \otimes e \leq x \otimes e \leq \top \otimes e = e\)
\end{example}

\begin{prop}[The “smooth” slice construction
    \cite{zé:good_quantales}]
    
    The slice construction above obviously does not preserve 
    unitality, integrality, (although it preserves 
    semicartesian-ness). When \(\quantale Q\) is 
    (left | right)-divisible and commutative (so that it is 
    immediately also integral) we can adjust \(\tensor\) on 
    each slice so as to remain divisible and commutative:
    If \(a,a'\in[0,b]\), we can define \(\tensor_b\) as follows
    \[
        b\tensor(b\residue a)\tensor(b\residue a')
    \]
\end{prop}
\begin{proof}
    It is not hard to see that in any quantale \(\quantale Q\),
    \[
        a\tensor(a\residue(a\tensor b)) = a\tensor b
    \]
    and in left-divisible quantales, one can then show that
    \[
        b\leq a \implies a\tensor(a\residue b) = b
    \]
    and similarly for right-divisible quantales (with the 
    appropriate residue). Since we are in a commutative setting 
    -- which is important later --, we can just assume the above 
    to hold.

    First, let us show that \(\tensor_b\) as defined is indeed 
    associative. We do this by deferring associativity to 
    \(\tensor\): 
    \begin{align*}
        x\tensor_b(y\tensor_b z) 
        &=  b\tensor(b\residue x)\tensor(b\residue[
            \overbrace{
                b\tensor(b\residue y)\tensor(b\residue z)
            }^{y\tensor_bz}
        ])\\
        &=  (b\residue x)\tensor b\tensor(b\residue[
            b\tensor(b\residue y)\tensor(b\residue z)
        ])\\
        &=  (b\residue x)\tensor{[
                b\tensor(b\residue y)\tensor(b\residue z)
        ]}\\
        &=  b\tensor
            (b\residue x)\tensor
            (b\residue y)\tensor
            (b\residue z)
    \end{align*}
    Now, we can repeat that but swapping \(x\) and \(z\); the 
    result is still the same -- of course --, yielding
    \[
        x\tensor_b(y\tensor_b z) = 
        z\tensor_b(y\tensor_b x) 
    \]
    But commutativity gives 
    \[
        z\tensor_b(y\tensor_b x) = 
        z\tensor_b(x\tensor_b y) = 
        (x\tensor_b y)\tensor_b  z
    \]
    Also, notice that \(\tensor_b\) is commutative. Now suppose 
    \(y, x_i\leq b\), then 
    \begin{align*}
        \left(\bigvee_{i} x_i\right)\tensor_by
        &=  b\tensor
            \left(
                b\residue\bigvee_i x_i
            \right)\tensor
            (b\residue y) \\
        &=  \left(
                \bigvee_i x_i
            \right)\tensor
            (b\residue y) \\
        &=  \bigvee_i x_i \tensor (b\residue y) \\
        &=  \bigvee_i b\tensor
                (b\residue x_i) \tensor 
                (b\residue y) \\
        &=  \bigvee_i x_i\tensor_b y
    \end{align*}
    which, together commutativity, then gives full 
    distributivity. Lastly, we should prove divisibility. 
    Suppose \(x \leq a\leq b\), so that we may try and find 
    a \(y\leq b\) such that \(x = y\tensor_b a\) -- and 
    commutativity gives full bidivisibility. Divisibility 
    on \(\quantale Q\) gives us that 
    \(a\tensor(a\residue x) = x\). Let's then consider 
    \[
        y = b \tensor (a \residue x)
    \]

    Bidivisibility gives \(y\leq b\) since it implies that 
    \(\quantale Q\) is semicartesian and obviously 
    \(b \wedge (a \residue x)\leq b\). Now, to show that the 
    tensor has the desired property:
    \begin{align*}
        a\tensor_b[b \tensor (a \residue x)]
        &=   b\tensor
            (b\residue a)\tensor
            (b\residue{[b\tensor (a \residue x)]})
        \\&=(b\residue a)\tensor
             b\tensor
            (b\residue{[b\tensor (a \residue x)]})
        \\&=(b\residue a)\tensor
            {[b\tensor (a \residue x)]}
        \\&= b\tensor
            (b\residue a)\tensor
            (a\residue x))
        \\&= a\tensor
            (a\residue x))
        \\&= x
    \end{align*}
\end{proof}
\begin{prop}
    Consider the category \(\Quantale_{c,d}\) of commutative 
    and divisible quantales and suprema, \(\top\) and 
    \(\tensor\) preserving morphisms between them. Now take a 
    one such quantale \(\quantale Q\). The smooth slice construction
    is, in fact, a functor 
    \[
        \quantale Q/\blank: \dual{\quantale Q}\to\Quantale_{c,d}
    \]
\end{prop}
\begin{proof}
    Take \(a\leq b\). First we need to provide a 
    morphism \(\quantale Q/b\to\quantale Q/a\) corresponding to the
    fact that \(a\leq b\). For that purpose, 
    \[
        (x\in\quantale Q/b)
        \xmapsto{\quantale Q/(a\leq b)}
        (a\tensor_b x\in\quantale Q/a)
    \]
    We know that \(b\mapsto a\tensor_b b = a\) and hence it 
    preserves \(\top\); we also know that 
    \[
        \bigvee_ix_i \mapsto a\tensor_b\bigvee_ix_i = 
        \bigvee_ia\tensor_bx_i
    \]
    and hence it preserves suprema. Finally, take 
    \(x\tensor_b y\), we know that it maps to 
    \[
        a\tensor_b(x\tensor_b y) = b\tensor
        (b\residue a)\tensor
        (b\residue x)\tensor
        (b\residue y)
    \]
    \begin{align*}
        (a\tensor_bx)\tensor_a(a\tensor_by)
        &=   a\tensor
            {[a\residue (a\tensor_bx)]}\tensor
            {[a\residue (a\tensor_by)]}
        \\&= a\tensor
            {[a\residue (b\tensor(b\residue a)\tensor(b\residue x))]}\tensor
            {[a\residue (a\tensor_by)]}
        \\&= 
            (b\tensor(b\residue a)\tensor(b\residue x))\tensor
            {[a\residue (a\tensor_by)]}
        \\&= 
            a\tensor(b\residue x)\tensor
            {[a\residue (a\tensor_by)]}
        \\&= 
            (b\residue x)\tensor(a\tensor_by)
        \\&= 
            (b\residue x)\tensor b\tensor
            (b\residue a)\tensor
            (b\residue y)
        \\&= b\tensor
            (b\residue x)\tensor 
            (b\residue a)\tensor
            (b\residue y)
    \end{align*}

    This proves it is a morphism, now it remains to be seen 
    that \(\quantale Q/\blank\) is actually functorial. The 
    first bit of functoriality is straightforward: 
    \(\quantale Q/(b\leq b)\) is \(b\tensor_b\blank\) which
    is the identity since \(b\) is the unit. Now take
    \(a\leq b\leq c\) 
    \begin{align*}
        [\quantale Q/(a\leq b)\circ
        \quantale Q/(b\leq c)](x) 
          &= 
        a\tensor_b(b\tensor_c x)
        \\&=    
         b\tensor
        (b\residue a)\tensor
        (b\residue[c\tensor(c\residue b)\tensor(c\residue x)])
        \\&=
         b\tensor
        (b\residue a)\tensor
        (b\residue{[b \tensor(c\residue x)]})
        \\&=
        (b\residue a)\tensor
        {[b \tensor(c\residue x)]}
        \\&=
        a\tensor(c\residue x)
        \\\\
        {[\quantale Q/(a\leq c)]}(x) 
        &=  a\tensor_c x
        \\&=c\tensor(c\residue a)\tensor(c\residue x)
        \\&=a\tensor(c\residue x)
    \end{align*}
\end{proof}

\begin{example}
    There is also a notion of product of quantales whose order
    and product are given point-wise. It is quite trivial to see 
    that the point-wise identity will be the identity, point-wise
    product will be the products etc.
\end{example}

\begin{remark}
    The last two examples introduced in 
    \ref{example: unital quantales} are neither commutative 
    nor semicartesian. The first is not idempotent but the 
    second is, and both are right-sided (resp. left-sided) 
    quantales 
    \cite[cf.]{rosenthal1990quantales}.
\end{remark}

\begin{example}
    The main examples of strong quantales are Heyting algebras, 
    \(([0,1],\leq,\cdot)\), strict linear quantales. 
    Some MV-Algebras, like the Chang's 
    \(([0,1],\wedge, \vee, \oplus, \tensor, 0, 1)\) \cite{chang}'s,
    are not strong.\\
\end{example}

\begin{remark}
    \begin{enumerate*}[label=(\roman*)]
        \item The class of commutative and semicartesian quantales 
            is closed under arbitrary products and and slices 
            \([a,b]\) whenever \(a\) is idempotent; 
        \item the class of strong quantales is closed under 
            arbitrary products and under interval constructions;
        \item as mentioned, commutative divisible quantales 
            are closed under smooth slices. 
    \end{enumerate*}
\end{remark}

\begin{center}
    \textbf{
        From now on, we assume all quantales to be commutative.
    }
\end{center}

\begin{definition}
    Let \(\quantale Q\) be a a quantale, we define an alternative 
    partial order \(\preceq\) given by
    \[
        a\preceq b \iff a = a \tensor b
    \]
\end{definition}

\begin{prop}[Properties for Semicartesian Quantales]
    \label{remark:semicartesian properties}
    Given a semicartesian quantale, 
    \begin{enumerate}
        \item 
            Let \(\quantale Q\) be a unital quantale, then it is 
            integral.
        \item 
            \(a\preceq b\leq c\) implies \(a\preceq c\);
        \item 
            If \(e \in \extent\quantale Q\), 
            \((e \preceq a) \iff (e \leq a)\). 
    \end{enumerate}
\end{prop}

\begin{prop}\label{prop:idempotent + semicartesian = locale}
    If \(\quantale Q\) is semicartesian and idempotent, it is in 
    fact a complete distributive lattice and \(\tensor=\wedge\). 
    In other words, it is a locale.
\end{prop}
\begin{proof}
    Suppose \(\quantale Q\) is in fact idempotent, we have --
    because \(\tensor\) is increasing in both arguments -- that
    \[
        a\leq b \implies a\leq c \implies a\leq b\tensor c
    \]

    Hence, if \(a\) is less than both \(b\) and \(c\), then it 
    must be smaller than \(b\tensor c\); but 
    since \(\quantale Q\) is semicartesian, by remark 
    \ref{remark:semicartesian properties} above, 
    \(\tensor\leq\wedge\). This means that \(b\tensor c\) is 
    a lowerbound for \(b\) and \(c\), but what we had gotten from 
    idempotency means it's the greatest such upper bound.

    Thus the multiplication satisfies the universal property of 
    infima. The above is just a particular case of 
    \cite[Proposition 2.1]{nlab:quantale}.
\end{proof}



\subsection{On \(\quantale Q\)-Sets} 

\begin{remark}
Hereon we are working exclusively with commutative    semicartesian quantales, as opposed to ones without those   properties.
\end{remark}

Given a quantale \(\quantale Q\), one may form -- roughly speaking -- a \(\quantale Q\)-space, wherein distances between points are measured by elements of \(\quantale Q\) as opposed to -- say -- \([0,\infty)\) as we often do. This definition is made precise in the notion of a \(\quantale Q\)-set.
\begin{definition}
    A \(\quantale Q\)-set is a set endowed with a    \(\quantale Q\)-distance operation usually denoted by    \(\delta\). \lat{i.e} a \(\quantale Q\)-set is a set \(X\) with  a map \(\delta:X^2\to X\) satisfying:
    \begin{enumerate}
        \item \(\delta(x,y) = \delta(y,x)\);
        \item \(\delta(x,y)\tensor\delta(y,z)\leq\delta(x,z)\);
        \item \(\delta(x,x)\tensor\delta(x,y) = \delta(x,y)\).
    \end{enumerate}
    and it is usual to denote \(\delta(x,x)\) by simply the   “extent of \(x\)” written as \(\extent x\).
\end{definition}

A couple of things might jump out to reader in the above definition. 
\begin{enumerate*}[label=(\roman*)]
    \item	\(\delta\) is symmetric, even though we have thrown out 
        all but the vaguest notions tying ourselves to metric 
        spaces;
    \item	Why is the triangle inequality upside down?
    \item	\(\extent x\tensor\delta(x,y) = \delta(x,y)\),
        why not just ask that \(\extent x = \top\)?
\end{enumerate*}

Those questions are all valid -- and answering the first and last ones 
differently has been done in the past and is the main difference 
between \(\quantale Q\)-sets and \(\quantale Q\)-enriched categories 
from a definitional perspective. The question of order being inverse is 
more one of sanity: since we treat a \(\quantale Q\)-set as a set with
\(\quantale Q\)-valued equality, it makes sense to think that 
\(\extent x\) is the maximally valid equality to \(x\) and hence the
triangular inequality needs to be turned upsidedown -- and turned into
the transitivity of equality.

\begin{remark}
    When we speak of properties of the type 
    \(P(\vec x)\leq Q(\vec x)\) in \(\quantale Q\)-sets, it is 
    often more insightful to think that the logically equivalent 
    (but notationally less helpful) statement 
    \[
        P(\vec x) \residue Q(\vec x) \ (=\top)
    \]
\end{remark}

There are two main category structures that one can canonically endow
the collection of all \(\quantale Q\)-sets with. One is taking maps to
be co-contractions (\lat{i.e.} they make \(\delta\) bigger) -- the 
other is to consider well behaved \(\quantale Q\)-valued relations 
between the underlying sets.

\begin{definition}
    A functional morphism \(f:X\to Y\) is a function \(f\) between
    the underlying sets of \(X\) and \(Y\) such that \(f\) 
    increases \(\delta\) and preserves \(\extent\); that is to say
    \[
        \delta_X \leq \delta_Y\circ (f\times f) 
    \]
    \[
        \extent_X = \extent_Y \circ f
    \]

    \begin{remark}
        There is a suitable notion of morphism between 
        \(\quantale Q\)-sets which is carried by 
        “\(\quantale Q\)-valued” relations. This notion isn't 
        explored in this paper, but they are called “relational
        morphisms”
    \end{remark}

    The reader should beware we don't often distinguish between 
    \(\delta_X\) and \(\delta_Y\) and instead rely on suggestively 
    named variables so as to indicate their type and hence the 
    \(\delta\) they refer to. In other words, the reader is 
    expected to be familiar with Koenig lookup%
    \footnote{
        ~Which, to quote a great website -- 
        \url{cppreference.com} -- 
        \begin{quote}
            Argument-dependent lookup, also known as ADL, 
            or Koenig lookup, is the set of rules for 
            looking up the unqualified function names in 
            function-call expressions, including implicit 
            function calls to overloaded operators.
        \end{quote}
    }.

    We denote by \(\QSet_r\) the category \(\quantale Q\)-sets 
    and relational morphisms between them and by \(\QSet_f\) the 
    category with the same objects but functional morphisms between
    them instead. Since we won't be tackling \(\QSet_r\), we take 
    \(\QSet\) to simply mean \(\QSet_f\).
\end{definition}

\begin{definition}
    Instead of proving the category axioms for functional morphisms
    we promised to prove a stronger result -- which is incidently 
    useful for another paper of ours in the works -- which is to 
    prove that \(e\)-morphisms form a category (given a generic
    commutative unital quantale) and that functional morphisms form
    a wide subcategory of \(I\)-morphisms (for \(I\) the monoidal
    unit).

    So, let us define \(e\)-morphisms: given an idempotent element 
    \(e\) of \(\quantale Q\), a \(e\)-morphism is a 
    functional morphism “up to error \(e\)”:
    \[
        e\tensor\delta(x,x') \leq \delta(f(x),f(x'))
    \]
    \[
        \extent f(x) = e\tensor\extent f(x)
    \]
\end{definition}
\begin{prop}
    We claim that the collection of \(\langle e,\phi\rangle\) where
    \(\phi\) is an \(e\)-morphism constitutes a category under the 
    obvious composition laws. Furthermore, the identity function is 
    a \(I\)-morphism where \(I\) is the unit of the quantale, and 
    further still: \(I\)-morphisms are closed under composition and
    form a subcategory which is definitionally equal to 
    \(\QSet_f\).
\end{prop}
\begin{proof}
    Firstly, the obvious composition takes an \(e\)-morphism \(f\)
    and an \(e'\)-morphism \(g\) to a \((e\tensor e')\) morphism
    \(g\circ f\). Associativity is due to functional (in \(\Set\), 
    that is) \(\circ\) associativity and the fact that \(\tensor\)
    makes \(\quantale Q\) a semigroup. The fact that \(g\circ f\) 
    is a \((e\tensor e')\)-morphism is rather obvious and the proof
    is ommited.
    
    The identity is evidently a \(I\)-morphism -- and of course 
    that composing \(I\)-morphisms gives a 
    \(I\tensor I = I\)-morphism.
\end{proof}

\subsection{Some Examples}

\begin{example} \label{empty-ex}
    The empty set is -- vacuously -- a \(\quantale Q\)-set.
\end{example}

\begin{example} \label{sep-ex}
    The set of idempotent elements of \(\quantale Q\) will be denoted 
    \(\extent\quantale Q\). Let \(X \subseteq \extent \mathscr{Q}\), then \((X,\otimes)\) is a \(\quantale Q\)-set. It is trivial to see 
    that \((e,e') \in \extent\quantale Q \times \extent\quantale Q \mapsto (e \otimes e') \) satisfies all \(\quantale Q\)-set laws, but as a first non-trivial example, we provide the details:
    \begin{align*}
    \delta(e,e')&= e \otimes e'\\
        &=e' \otimes e\\
        &=\delta(e',e)\\
        \\
        \delta(e,e')\otimes \delta(e',e'')&= e \otimes e' \otimes e' \otimes e''\\
        &\leq e \otimes e''\\
        &= \delta(e',e)\\
        \\
        \delta(e,e')\otimes \extent e = & e \otimes e' \otimes e \otimes e\\
        &= e \otimes e'\\
        &= \delta(e,e')\\
    \end{align*}
    Note that  for any \(e \in X\)  satisfies \(\delta(e,e) = \extent e =  e\) and that   \(e \otimes e'\) is the infimum of \(\{e,e'\}\) in the poset \( \extent\quantale Q \).
\end{example}

\begin{remark}
    One cannot use the whole of \(\quantale Q\) in 
    \(\extent\quantale Q\)'s stead, as 
    \(\delta(x,y)\tensor\extent x\) would not hold. The only reason 
    it holds in the above example is because for idempotents 
    \(\tensor=\wedge\). However, under conditions, one can obtain a \(\quantale Q\)-set
    where the underlying set is \(\quantale Q\) itself.
\end{remark}

\begin{example} \label{Qsets-ex}
    Suppose that \(\quantale Q\) is a quantale with 
    “idempotent upper  approximations”\footnote{%
        In \cite{tenorio2022sheaves} are described sufficient conditions
        for \(\quantale Q\) to have such a property.
    }: 
    \[
        \forall q\in\quantale Q:\exists q^+\in\extent\quantale Q:
            q^+ = \min\set{e\in\extent\quantale Q}{q\preceq e}
    \]
    Then
    \[
        \delta(x, y) = \begin{cases}
            x\tensor y, & x\not=y;\\
            x^+,        & x  =  y.
        \end{cases}
    \]
    defines a \(\quantale Q\)-set structure on \(\quantale Q\) itself.
\end{example}

\begin{example}
    Much akin to how Lawvere's quantale is a Lawvere space, a 
    (integral and commutative) quantale \(\quantale Q\) is a 
    \(\quantale Q\)-set. This is achieved with the following:
    \[
        \delta(x,y) = (x\residue y)\wedge(y\residue x)
    \]
    which is roughly equivalent to \(|x - y|\) for real numbers.
    This isn't necessarily the best \(\quantale Q\)-set structure
    we can give them, as \(\extent x = \top\) for any \(x\). 

    Ways to mitigate this phenomenon, which is specially strange 
    for \(\bot\), involve taking into account idempotents above 
    \(x\). An important quantallic property is the existence of 
    an operation \((\blank)^-\) taking an element \(x\) to the value 
    \(\sup\set{e\in\extent\quantale Q}{e\preceq x}\). Multiplying 
    \(\delta(x,y)\) by \(x^-\tensor y^-\) guarantees -- for instance
    -- that the above construction coincides with the terminal 
    object when \(\quantale Q\) is a locale.

    Another way to correct this, is to incorporate 
    \(\extent\quantale Q\) more directly, considering the space
    with underlying set \(\quantale Q \times\extent\quantale Q\)
    and \(\delta\) given by
    \[
        \delta((x,e),(y,a)) = 
        a\tensor e\tensor{[(x\residue y)\wedge(y\residue x)]}
    \]
    We write this \(\quantale Q\)-set as \({\quantale Q}_{\extent}\).
\end{example}

\begin{example}
    A construction that is explored in this work 
    \cite{AMM2023constructions} but is suited to be mentioned here
    is \(X\chaosor X\), given by the underlying set 
    \(|X|\times|X|\) and with \(\delta\) given by the product of 
    the \(\delta\) of the coordinates. The reason this construction
    is relevant here is because 
    \(\delta: X\chaosor X\to{\quantale Q}_{\extent}\)
    in a natural way:
    \[
        (x,y)\mapsto (\delta(x,y), \extent x\tensor\extent y)
    \]
    And this happens to be a functional morphism.
\end{example}

\begin{example}
     Suppose \((X,d)\) is a pseudo-metric space, 
     then \((X,d)\) is a \([0,\infty]\)-set where 
     \(\delta(x,y) \neq \bot = \infty, \forall x,y \in X\).
\end{example}


\begin{example}
    Given a commutative ring \(A\), let the set of its (left) ideals 
    be denoted \newcommand\ideals{\mathscr I} \(\ideals_A\). 
    \(\ideals_A\) is a quantale. Given a left \(A\)-module \(M\), we 
    can endow it with the structure of a \(\ideals_A\)-set:
    \[
        \delta(x,y) = \bigvee\set{
            I \in\ideals_{A}
        }{
            I\cdot x = I\cdot y
        }
    \]
    In fact, that supremum is attained at with a particular ideal. 
    Moreover, \(Ex = A = \max\ideals_A\). 
\end{example}


\begin{remark}[Completeness]
    As mentioned previously, there is a category of \(\quantale Q\)-sets 
    with relational morphisms between them. It happens that this category 
    is (equivalent to) a reflective subcategory of \(\QSet\). The objects 
    of this reflective subcategory are called “Scott complete” 
    \(\quantale Q\)-sets.

    The notion of completeness has to do with singletons, which are 
    \(\quantale Q\)-valued distributions which “measure a point”. And 
    being Scott-complete the same as saying that all points we measure 
    are actually there, and we don't measure any points more than once.

    This, in a sense, is a similar condition to that of a space being 
    sober. There is at least one different notion of completeness, called
    gluing completeness -- which is related to compatible local data 
    having exactly one gluing. These notions, and the reflective 
    subcategories they define are explored in a different article in 
    development \cite{AMM2023completeness}.
\end{remark}
	\section{Main Constructions}

In the present section, we provide the basic constructions in the category of all  \(\quantale Q\)-sets and functional morphisms: limits, colimits, (regular) suboject classifier, etc.

We start with the following:

\begin{prop}[Separating Family]
    For each \(e \in \extent\quantale Q\)   we 
    can consider the singleton set \(S_e = \{e\}\) endowed with the natural 
    \(\quantale Q\)-set structure described in Example\ref{sep-ex}.
    
    Then the set of such \(\quantale Q\)-sets is a separating family
    for the category of \(\quantale Q\)-sets.
\end{prop}
\begin{proof}
    Indeed, for any pair of morphisms \(f, g:X \to Y\), 
    if \(f \not= g\) then must be some \(x \in X\) such that 
    \(f(x) \not= g(x)\). Taking \(e = \{\extent x\} \in  \extent\quantale Q\) and the function \(s_x : S_e \to X \) given by \((e\in\{e\})\mapsto(x\in X)\), then  it is a functional morphism and  obviously separates 
    \(f\) from \(g\).
\end{proof}

\subsection{Limits}

\begin{prop}[Terminal Object]
    The terminal object is \(\top=(\extent \mathscr{Q},\otimes)\).
\end{prop}
\begin{proof}
    The set of idempotent elements of \(\quantale Q\), 
    \(\extent\quantale Q\), endowed with the 
natural 
    \(\quantale Q\)-set structure described in Example\ref{sep-ex} ( \(\delta(e,e') = e \otimes e'\) ) must be the 
    terminal object because \(\extent e = e\). 

For each   \(\quantale Q\)-set \((X, \delta)\), the unique morphism with codomain \(\top\) is defined as:
    \begin{align*}
        !:X &\to\top\\
        f(x)&=\extent x
    \end{align*}
    
    Since functional 
    morphisms preserve extents, one has that 
    \(\extent f(x) = \extent x\) -- however, \(\extent f(x) = f(x)\) 
    and thus \(f(x) = \extent x\). This proves that there is at most
    one morphism \(X\to\extent\quantale Q\).

    On the other hand, \(\delta(x,y)\leq\extent x\otimes\extent y\) 
    which just happens to make \(x\mapsto\extent x\) a functional 
    morphism:

     \begin{align*}
        \delta_X(x,x')\leq \extent x &\otimes \extent x' = \delta(x,x')
        \\
        \extent x = \extent x &\otimes \extent x = \extent x
    \end{align*}
 
\end{proof}

\begin{prop}[Non-empty products]
    The product of \(\quantale Q\)-sets \((X_i)_{i \in I}\) is:
\end{prop}
\begin{proof}
    \[
        \prod_{i \in I} X_i
            = \left(\set{
                (x_i)_{i \in I} \in \prod_{i \in I} |X_i|
            }{
                \forall i, j \in I (\extent x_i = \extent x_j)
            }, \delta
        \right)
    \]\[
        \delta((x_i)_{i \in I},(y_i)_{i \in I})
            = \bigwedge_{i \in I} \delta_i(x_i,y_i)
    \]
    \begin{remark}
        Since
        \[
            (x_i)_{i \in I} \in \prod_{i \in I} X_i
            \implies 
            \forall i, j \in I: \extent x_i = \extent x_j
        \]
        we can conclude that 
        \[
            \extent(x_i)_{i \in I} 
            = \bigwedge_{i \in I} 
                \extent x_i 
            = \bigwedge_i \extent x_i
        \]
    \end{remark}
\end{proof}
\begin{proof}[Proof: \(\quantale Q\)-set]
    \begin{align*}
                \delta((x_i)_{i \in I},(y_i)_{i \in I})
              &=\bigwedge_{i \in I} \delta_i(x_i,y_i)
            \\&=\bigwedge_{i \in I} \delta_i(x_i,y_i)
            \\&=\delta((y_i)_{i \in I},(x_i)_{i \in I})
        \\\\
                \delta((x_i)_{i \in I},(y_i)_{i \in I})
                \otimes
                \delta((x_i)_{i \in I},(x_i)_{i \in I})
            \\&=\bigwedge_{i \in I}\delta(x_i,y_i)
                \otimes
                \bigwedge_{j \in I}\delta(x_j,x_j)
            \\&=\bigwedge_{i \in I}\delta(x_i,y_i) 
                \otimes
                \extent(x_j)_{j \in J}
            \\&=\bigwedge_{i \in I}\delta(x_i,y_i)\otimes\extent(x_i)
            \\&=\bigwedge_{i \in I}\delta(x_i,y_i)
        \\\\
            \delta((x_i)_{i \in I},(y_i)_{i \in I})
            \otimes
            \delta((y_i)_{j \in I},(x_i)_{z \in I})
          &=\bigwedge_{i \in I}
                \delta(x_i,y_i)
                \otimes
                \bigwedge_{j \in I}\delta(y_j,z_j)
        \\&\leq\bigwedge_{i \in I}
            \delta(x_i,y_i)\otimes\delta(y_i,z_i)
        \\&\leq\bigwedge_{i \in I}
            \delta(x_i,z_i)
        \\&=\delta((x_i)_{i \in I},(y_i)_{i \in I})
    \end{align*}
\end{proof}
\begin{proof}[Proof: Projections]
    \begin{align*}
        \pi_i:\prod_{i \in I} X_i \to& X_i\\
        \pi_i((x_j)_{j \in I})&=x_i
    \end{align*}
    which are morphisms because
    \begin{align*}\
        \delta_i(\pi_i((x_j)_{j \in I}),\pi_i((y_j)_{j \in I}))
            &=      \delta_i(x_i,y_i)\\
        \\  &\geq   \bigwedge \delta_j(x_j,y_j)\\
        \\  &=      \delta((x_j)_{j \in I},(y_j)_{j \in I})
        \\\\
        \delta_i(\pi_i((x_j)_{j \in I}),\pi_i((x_j)_{j \in I}))
            &= Ex_i
        \\  &= \bigwedge_{j \in I} Ex_j
        \\  &= \delta((x_j)_{j \in I},(x_j)_{j \in I})
    \end{align*}
\end{proof}
\begin{proof}[Proof: Universality]
    Let \(f_i:A\to X_i\) a family of morphisms. We define:
    \begin{align*}
        h:A &\to \prod_{i \in I} X_i\\
        h(a)&=(f_i(a))_{i \in I}
    \end{align*}
    Then \(h\) satisfies the universal property:
    \begin{align*}
        \pi_i\circ h(a) &= \pi_i ((f_i(a))_{i \in I}) = f_i(a)\\
        \pi_i\circ h    &= f_i
    \end{align*}
\end{proof}

\begin{prop}
    The equalizer of \(f,g:(X,\delta) \to (Y,\delta')\) is:
    \[
        \equalizer (f,g) = \left(
            \set{
                x \in X
            }{
                f(x)=g(x)
            },
            \delta\restriction_{\set{x \in X}{f(x)=g(x)}}
        \right)
    \]
\end{prop}
\begin{proof}
    It is obviously a \(\quantale Q\)-set. It obviously equalizes the pair,
    now, for a given \(\alpha:A\to X\) also equalizing the pair, 
    Realize that \(\alpha\)'s image must be entirely contained in 
    \(\equalizer(f,g)\) as 
    \[
        [f\circ\alpha(a) = g\circ\alpha(a)]\iff a\in\equalizer(f,g)
    \]
    and hence the restriction of codomain yields the unique morphism that
    makes the diagram commute and establishes universality.
\end{proof}

\begin{remark}
    Considering that limits are given as equalizers of products, the 
    category has been shown to be complete.
\end{remark}

\begin{prop} [Monomorphisms]
    A morphism is a monomorphism iff it is a injective as a function.
\end{prop}
\begin{proof}
    Injectivity implying monomorphicity is trivial. The only reason the 
    other direction isn't trivial is that, in principle, there might not 
    be enough morphisms to witness a non-injective morphism not being a 
    mono. But this is easy: suppose \(f(x) = f(y)\), consider the 
    \(\quantale Q\)-set \(\{*_e\}\) for \(e = \extent x = \extent y\), 
    with \(\extent *_e = e\). This \(\quantale Q\)-set (eg. \ref{sep-ex})
    can be included in our object by sending \(*_e\) to either \(x\) or 
    \(y\); if \(f\) is monic, then both we would have that \(x = y\) 
    thus establishing injectivity.
\end{proof} 

\begin{prop} [Regular monomorphisms]
   The class of regular monomorphisms coincides with the class of 
   injective morphisms that preserve \(\delta\). Which are precisely
   the \(\delta\)-preserving injective functions between the underlying
   sets.
\end{prop}
\begin{proof}
    A monomorphism is regular when it is an equalizer of a pair of 
    parallel arrows. Suppose \(g,h:A\to B\), one way to conceive of 
    an equalizer is as being the maximal subobject of \(A\) making 
    the diagram commute. It is quite trivial to see that subobjects
    in general \emph{must} be subsets with a point-wise smaller 
    \(\delta\). Hence, the largest subobject of \(A\) that 
    equalizes the pair is simply the subset of \(A\) where the 
    functions agree on, with the largest possible \(\delta\): that 
    being \(A\)'s. 

    Hence, we see a pattern where equalizers -- in general -- are 
    simply subsets with \(\delta\) coming from a restriction of 
    the original object. Importantly, though, we refer to 
    “regular subobjects” as monomorphisms that preserve \(\delta\)
    -- as they have been equivalently characterized.

    The skeptical reader might question that we have merely shown
    that regular monos preserve \(\delta\), as opposed to showing 
    such to be a sufficient condition. 
    
    In which case, given the fact that monos are injective 
    functions, \(\delta\)-preserving monos are simply subsets with 
    its superset's \(\delta\); consider one such mono: 
    \(f :A\rightarrowtail X\) (we may think that
    \(A\subseteq X\) and \(\delta_A = \delta_{X_{| A \times A}}\)), one takes 
    \(X_f = (X\amalg X)/\sim_f\) with \(\sim_f\) defined so as to exactly identify 
    both copies of \(A\). 
    \[
        \delta(\eqclass{(x,i)},\eqclass{(y,j)}) = 
        \begin{cases}
            \delta(x,y), & i = j\\
            \bigvee_{a\in A}\delta(x,f(a))\tensor\delta(f(a),y), 
                & i\not = j
        \end{cases}
    \]
\end{proof}
\begin{proof}[Proof: \(X_f\) is a \(\quantale Q\)-set.]
    Take \(i,j,k\) to be different indices
    \begin{align*}
        \delta(\eqclass{(y,i)},\eqclass{(y',i)})&=\delta(y,y')\\
        &=\delta(y',y)\\
        &=\delta(\eqclass{(y',i)},\eqclass{(y,i)})\\
    \end{align*}
    \begin{align*}
        \delta(\eqclass{(y,i)},\eqclass{(y',i)}) \otimes \delta(\eqclass{(y',i)},\eqclass{(y'',i)})&=\delta(y,y') \otimes \delta(y',y'')\\
        &\leq\delta(y,y'')\\
        &=\delta(\eqclass{(y,i)},\eqclass{(y'',i)})\\
    \end{align*}
    \begin{align*}
        \delta(\eqclass{(y,i)},\eqclass{(y',i)}) \otimes \delta(\eqclass{(y,i)},\eqclass{(y,i)})&=\delta(y,y')\otimes \delta(y,y)\\
        &=\delta(y,y')\\
        &=\delta(\eqclass{(y,i)},\eqclass{(y',i)})\\
    \end{align*}
    \begin{align*}
        \delta(\eqclass{(y,i)},\eqclass{(y',j)})&=\bigvee_{x \in X}\delta(y,f(x))\otimes \delta(f(x),y')\\
        &=\bigvee_{x \in X}\delta(y',f(x))\otimes \delta(f(x),y)\\\\
        &=\delta(\eqclass{(y',j)},\eqclass{(y,i)})\\
    \end{align*}
    \begin{align*}
        \delta(\eqclass{(y,i)},\eqclass{(y',j)}) \otimes \delta(\eqclass{(y',j)},\eqclass{(y'',j)]}&=\bigvee_{x \in X} \delta(y,f(x)) \otimes \delta(f(x),y') \otimes \delta(y',y'')\\
        &\leq \bigvee_{x \in X} \delta(y,f(x)) \otimes \delta(f(x),y'1)\\
        &=\delta(\eqclass{(y,i)},\eqclass{(y'',j)})\\
    \end{align*}
    \begin{align*}
        \delta(\eqclass{(y,i)},\eqclass{(y',j)}) \otimes \delta(\eqclass{(y',j)},\eqclass{(y'',k)})&=\bigvee_{x \in X}\bigvee_{x' \in X} \delta(y,f(x)) \otimes \delta(f(x),y') \otimes \delta(y',f(x'))\otimes \delta(f(x'),y'')\\
        &\leq\bigvee_{x \in X}\bigvee_{x' \in X} \delta(y,f(x)) \otimes \delta(f(x),f(x')) \otimes\delta(f(x'),y'')\\
        &=\bigvee_{x \in X}\delta(y,f(x)) \otimes\delta(f(x),y'')\\
        &=\delta(\eqclass{(y,i)},\eqclass{(y'',k)})\\
    \end{align*}
\end{proof}
\begin{proof}[Proof: Regularity]
    Define \(g_0,g_1:Y \to X_f\) as \(g_i(y)=(y,i)\). 
    Both are morphisms:
    \[
        \delta(g_i(y),g_i(y'))=\delta((y,i),(y',i))=\delta(y,y')
    \]
    The equalizer of \(g_0\) and \(g_1\) are precisely \(\img f\):
    \begin{align*}
        g_0(y) = g_1(y) 
        &\implies \eqclass{(y,0)}=\eqclass{(y,1)}\\
        &\implies y \in Im f
    \end{align*}
    Thus, \(f\) is regular.
\end{proof}

\subsection{Colimits}

\begin{prop}[Initial object]
    The initial object is \(\bot = (\varnothing,\varnothing)\)
\end{prop}
\begin{proof}
    This is a \(\quantale Q\)-set by vacuity (Example \ref{empty-ex}). The empty function is a morphism
    again by vacuity, and it obviously must be the initial object, as morphisms
    are -- in particular -- functions as well.
\end{proof}

\begin{prop}[Non-empty coproducts]
    The coproduct of \(\quantale Q\)-sets \((X_i)_{i \in I}\) is:
    \begin{align*}
        \coprod_{i \in I} X_i&=(\coprod_{i \in I}|X_i|, \delta)\\
        \delta((x,i),(y,i))  &= \begin{cases}
            \delta_i(x,y), &i = j\\
            \bot,          &i\not= j. 
        \end{cases}
    \end{align*}
\end{prop}
\begin{proof}[Proof: It is a \(\quantale Q\)-set]
    Suppose that \(i\neq j \neq k \neq i\). Without loss of generality:
    \begin{align*}
        \tag{\(\delta\) is symmetric}
        \delta((x,i),(y,i))&=\delta_i(x,y)\\
        &=\delta_i(y,x)\\
        &=\delta((y,i),(x,i))\\
        \\
        \delta((x,i),(y,j))&=\bot\\
        &=\delta((y,j),(x,i))\\
        \\
        \tag{extensionality}
        \delta((x,i),(y,i))\tensor E(x,i) &= \delta_i(x,y) \tensor Ex\\
        &= \delta_i(x,y)\\
        &= \delta((x,i),(y,i))\\
        \\
        \delta((x,i),(y,j)\tensor E(x,i) &= \bot\tensor Ex\\
        &= \bot\\
        &= \delta((x,i),(y,j))\\
        \\
        \tag{triangular inequality}
        \delta((x,i),(y,i))\tensor \delta((y,i),(z,i)) &= \delta_i(x,y)\tensor \delta_i(y,z)\\
        &\leq\delta_i(x,z)\\
        &=\delta((x,i),(z,i))\\
        \\
        \delta((x,i),(y,j)\tensor\delta((y,j),(z,j))&=\bot\tensor\delta((y,j),(z,j)\\
        &= \bot\\
        &\leq \delta((x,i),(z,j))\\
        \\
        \delta((x,i),(y,j)\tensor \delta((y,j),(z,k)) &= \bot\tensor \bot\\
        &= \bot\\
        &\leq \delta((x,i),(z,k))
    \end{align*}
\end{proof}
\begin{proof}[Proof: Coprojections]
    \begin{align*}
        \coproj_i: X_i\to& \coprod_{i \in I} X_i \\
        \coproj_i(x_i)   &= (x_i,i)
    \end{align*}
    By construction those are obviously morphims.
\end{proof}
\begin{proof}[Proof: Universality]
    Let \(f_i:X_i\to A\) a family of morphisms. We define:
    \begin{align*}
        h: \coprod_{i \in I} X_i &\to A\\
        h(x,i)                   & =  f_i(x)
    \end{align*}
    Then \(h\) satisfies the universal property:
    \begin{align*}
        f_i\circ\coproj_i(x,i)=f_i(x)&=h(x,i)\\
        f_i\circ\coproj_i&=h
    \end{align*}
\end{proof}

\begin{example}
    Given a nonempty index set \(I\), we have a 
    \(\quantale Q\)-set \(\coprod_{i\in I}\top\), given by
    \[
        \delta((e, i), (e', i')) = \begin{cases}
            e\otimes e', &\text{if} ~ i = i';\\
            \bot,       &\text{otherwise}.
        \end{cases}
    \]    
\end{example}

\begin{prop}[Coequalizers]
    For $f,g:(X,\delta) \to (Y,\delta')$, we define the equivalence relation $\sim \ \subseteq |Y|\times |Y|$ as the transitive closure of $f(x)\sim f(x)\sim g(x)\sim g(x)$.

    The coequalizer of $f,g:(X,\delta) \to (Y,\delta')$ is:\\
    \begin{align*}
    \coequalizer (f,g)&=\left(\faktor Y \sim,\delta\right)\\
    \delta(\eqclass{y},\eqclass{y'})&= \bigvee_{\substack{a\sim y\\a'\sim y'}}      
        \delta(a,a') 
    \end{align*}
\end{prop}
\begin{proof}[Proof: \(\quantale Q\)-set]
    \begin{align*}
        \delta(\eqclass{y},\eqclass{y'})
        &=  \bigvee_{\substack{a\sim y\\a'\sim y'}} 
                \delta(a,a')
        =   \bigvee_{\substack{a\sim y\\a'\sim y'}} 
            \delta(a',a)
        \\&=\delta(\eqclass{y'},\eqclass{y})
    \end{align*}
    \begin{align*}
        \delta(\eqclass{y},\eqclass{y'})
        &\geq
            \delta(\eqclass{y},\eqclass{y'})
            \tensor
            \extent\eqclass{y}
        \\&=\delta(\eqclass{y},\eqclass{y'})
            \tensor 
            \bigvee_{\substack{a\sim y\\a'\sim y}}\delta(a',a)
        \\&=\bigvee_{\substack{a\sim y\\a'\sim y'}}\delta(a',a)
                \tensor 
            \bigvee_{a\sim y} Ea 
        \\&\geq
            \bigvee_{\substack{a\sim y\\a'\sim y'}} 
                \delta(a',a)\tensor Ea
        \\&=\bigvee_{\substack{a\sim y\\a'\sim y'}} \delta(a',a)
        \\&=\delta(\eqclass{y},\eqclass{y'})
    \end{align*}
    \begin{align*}
        \delta(\eqclass{y},\eqclass{y'})
        \tensor 
        \delta(\eqclass{y'},\eqclass{y''})
        &=  \bigvee_{\substack{a\sim y\\a'\sim y'}} \delta(a',a) \tensor
            \bigvee_{\substack{b\sim y'\\b'\sim y''}} \delta(b,b')
        \\&=\bigvee_{\substack{a\sim y\\a'\sim y'}}
            \bigvee_{\substack{b\sim y'\\b'\sim y''}}
                \delta(a',a) \tensor \delta(b,b')
        \\&=\bigvee_{\substack{a\sim y\\a'\sim y'}}
            \bigvee_{\substack{b\sim y'\\b'\sim y''}}
                \delta(a',a) \tensor Ea \tensor \delta(b,b')
        \\&=\bigvee_{\substack{a\sim y\\a'\sim y'}}
            \bigvee_{\substack{b\sim y'\\b'\sim y''}}
                \delta(a',a) \tensor 
                \delta(a,b)  \tensor
                \delta(b,b')
        \\&\leq
            \bigvee_{\substack{a\sim y\\a'\sim y'}}
            \bigvee_{\substack{b\sim y'\\b'\sim y''}} 
                \delta(a,b) \tensor \delta(b,b')
        \\&\leq
            \bigvee_{\substack{a\sim y\\a'\sim y'}}
            \bigvee_{\substack{b\sim y'\\b'\sim y''}} 
                \delta(a,b')
        \\&=\bigvee_{\substack{a\sim y\\b'\sim y'}} 
            \delta(a,b')
        \\&=\delta(\eqclass{y},\eqclass{y''})
    \end{align*}
\end{proof}

\begin{prop}[Epimorphisms]
    Epimorphisms are precisely the surjective morphisms. 
    Suppose that $f:X \to Y$ is surjective. Then, 
    $\forall y \in Y: \exists x \in X: f(x)=y$. Let $g,h:Y\to A$. 
    Suppose that $\forall x \in X :g \circ f(x)=h \circ f(x)$. 
    Let $y \in Y$. There is a $x \in X$ such that $f(x)=y$. 
    Then $g(y)=g \circ f(x)=h\circ f(x)=h(y)$.\\
\end{prop}
\begin{proof}
    Take some \(f:X\to Y\) and suppose
    \[
        \exists y_0\in Y:\forall x\in X: f(x)\not = y_0
    \]
    Let \(Y_f=\left(\faktor{(Y\amalg Y)}\sim, \delta\right)\), 
    as generated by 
    \[
        y\in\img f \implies \coproj_l(y) \sim \coproj_r(y)
    \]
    Let $g,h:Y\to Y_f$ as $g(y)=\eqclass{(y,0)}$ and 
    $h(y)=\eqclass{(y,1)}$. 
    Then $g\neq h$ but $g \circ f = h \circ f$.
\end{proof}

\subsection{Subobject Classifier}

The goal of this short subsection is to establish that the category \(\quantale Q\)-sets has a (almost trivial) classifier for the regular subobjects.

\begin{remark} 
    Note that the category  is not balanced, since there are many bijective morphisms -in particular, morphisms that are mono+epi- that are not isomorphisms, that coincides with the bijective morphisms that preserves \(\delta\). 
\end{remark}

\begin{prop}[Regular Subobject Classifier]
    Let \(\Omega=(\top\mathbin{\dot\cup}\top,\delta)\) where:
    \[
        \delta((e,i),(e',j))= e\otimes e'
    \]
    and consider the morphism \(t: \top \to \Omega\) that includes \(\top\) in the second copy of \(\top\) in \(\Omega\):  \(t(e) = (e,1)\). Then \(t: \top \to \Omega\) is a classifier for the regular subobjects.
\end{prop}
\begin{proof}
    Note first that \(\Omega\) is a well-defined \(\quantale Q\)-set and that the identity map \(\top\coprod\top \to \Omega\) is a bijective morphism that almost never is an isomorphims (is isomorphism iff \(\quantale Q\) has a unique idempotent member). Moreover, note that \(t: \top \to \Omega\) is a regular monomorphism. 
    
    For each regular monomorphism \(f:X \to Y\), we define \(\chi_f:Y\to \Omega\) as 
    \[
        \begin{cases}
            (\extent y,1), & y = f(x), \\
            (\extent y,0), & y \in Y \setminus f[X]
        \end{cases}
    \]
    It is evident that this is a morphism, as this is akin to the terminal arrow, but we plug in an extra tag that doesn't interfere with \(\delta\) but allows us to keep track of the element's provenance.

    \begin{enumerate}[label={\bf Claim}]
        \item $\chi_f \circ f = t \circ {!_X}$:\\
            \[
                \chi_f (f(x)) = (E f(x), 1) = (E x, 1) = t (Ex) = t (!_X(x))
            \]
        \item The commutative diagram above
        $$
            (\Omega \overset{t}\leftarrow \top \overset{!_X}\leftarrow X \overset{f}\to Y \overset{\check{f}}\to \Omega)
        $$ 
        is a pullback square.

        Let $u : X \to \top \underset\Omega\times Y$  be the unique morphism given by the universal property of pullbacks. We will show that $u$ is a bijective morphism that preserves $\delta$s, thus it is an isomorphism. 
        
        Note that  $u(x) = (E_X x, f(x))$ for each $x \in X$, since $E_\top (E_X x) = E_X x = E_Y f(x)$ and $t(E_X x) = \chi_f (f(x))$.

        $u$ is injective: If $u(x) = u(x')$ then, $x = x'$ since $f$ is injective. 
        
        $u$ is surjective: if $(e, y) \in \top \underset\Omega\times Y $, then $E_\top (e) = e = E_Y y$ and $(e,1) = t(e) = \chi_f (y)$ thus, by the definition of $\chi_f$, $y = f(x)$ for some (unique) $x \in X$. Then $e = E_Y y = E_Y f(x) = E_X x$ and $(e, y) = (E_X x, f(x))$. 

        $u$ preserves $\delta$s: 
        $$
            \delta(u(x), u(x')) = \delta ((E_X x, f(x)), (E_X x', f(x'))) =
        $$
        $$
            \delta_\top (E_X x, E_X x') \wedge \delta_Y (f(x), f(x')) = \delta_\top (E_X x, E_X x') \wedge \delta_X (x, x') =
        $$
        $$ 
            (E_X x \tensor E_X x') \wedge \delta_X (x, x') = \delta_X (x, x')
        $$ 
        finishing the proof of the claim.        

        \item $\chi_f$ is the unique arrow $\check{f} : Y \to \Omega$ such that the diagram 
        $$(\Omega \overset{t}\leftarrow \top \overset{!_X}\leftarrow X \overset{f}\to Y \overset{\check{f}}\to \Omega)  $$ 
        is a pullback.

        Let $x \in X$. By the commutativity of the diagram, $\check{f}(f(x)) = t(!_X(x)) = (E_X x,1) = (E_Y f(x),1)$. 

        Let $y \in Y \setminus f[X]$. Consider the $\quantale Q$-set $Z = \{e\}$, where $e = E_Y y$, as described in \ref{sep-ex} and the (well defined) morphism $g : Z \to Y$ such that $g(E_Y y) = y$. If  $\check{f}(g(e)) = t (!_Z(e)) = t (e) = (e,1)$, then there is a unique morphism $v : Z \to X$ such that $f(v(e)) = g(e) = y$. Thus $\check{f}(y) = \check{f}(g(e)) \neq (e,1)$. But $E_\Omega \check{f}(y) = E_\Omega \check{f}\circ g(e) = E_Z e = e $ and $\check{f}(y) \in  \{(e',0), (e',1)\}$ for some idempotent $e'$ such that $e' = e' \tensor e' = E_\Omega \check{f}(y) $. This means that $\check{f}(y) = (e, 0) = (E_Y y, 0)$, finishing the proof of the claim.
        
    \end{enumerate}

\end{proof}

	\section{Local Presentability}

Needless to say that a category being locally presentable is a very
strong property, in that it -- for instance -- allows us to construct
right adjuncts to any functor that “ought” to have them 
(cocontinous) when the categories in question are locally presentable.
It also reflects positively into the slices of the category in 
question.

Consider a setting in which we have an object \(X\) of \(\category C\)
-- a sufficiently cocomplete category -- and a diagram 
\(D:\category X\to\category C\). There is a canonical map, given by 
universality of the coprojections%
\footnote{
    Which the reader may read as “ip”, as it is the dual of pi.
} \(\coproj_k:D_k\to\colim_k D_k\):
\[\begin{tikzcd}
	{\colim\limits_k\hom(X,D_k)} & {\hom(X,\colim\limits_kD_k)} \\
	{\eqclass{(f, k)}} & {\coproj_k\circ f}
	\arrow["\phi", from=1-1, to=1-2]
	\arrow[maps to, from=2-1, to=2-2]
\end{tikzcd}\]

For \(\hom(X,\blank)\) to preserve colimits it then suffices that 
this natural map be a natural isomorphism. And hence, simply an 
isomorphism pointwise -- as it is already natural. This is how we 
shall go about showing that \(\QSet\) is accessible, and being 
accessible and cocomplete it must be locally presentable.

For accessibility, we need a regular cardinal \(\kappa\) -- to be 
determined -- and an (essentially) small collection of \(\kappa\)
compact objects which generate \(\QSet\) under \(\kappa\) directed
colimits. First, then, we are going to search for one such class of 
objects, then show they generate the category appropriately.

\subsection{\(\kappa\)-Compact Objects} 

To determine them, we must obviously settle on \emph{some} regular 
cardinal. It turns out that \(\left|\quantale Q\right|^+\) -- the 
successor cardinal of \(\quantale Q\)'s cardinality -- suffices. 
Thus defined, it is reasonably straightforward to present an 
essentially small class of \(\kappa\)-compact objects. 

Consider \(X\) a \(\quantale Q\)-set such that its carrier set has 
cardinality less than \(\kappa\). It -- we shall show -- is compact
wrt. \(\kappa\). This, of course, is to show that \(\hom(X,\blank)\)
preserves \(\kappa\)-directed colimits. This we do in two steps:
showing surjectivity and injectivity of \(\phi\) for our \(X\).

\begin{lemma}[\(\phi\) is surjective]
    Suppose \(D:\langle I,\leq\rangle\to\QSet\) is a
    \(\kappa\)-directed diagram -- or it could just be 
    \(\kappa\)-filtered it really doesn't matter either way.

    To show our claim, it would suffice that any arrow 
    \(X\to\colim_k D_k\) actually factors through some 
    (dependent on it) \(D_i\) -- as we can then take the canonical
    maps from those \(\hom\)s into \(\colim_k\hom(X,D_k)\) and have
    a section for \(\phi\) -- showing it is surjective.

    We do this by constructing an index set \(J\) such that a given 
    \(f\) must factor through any \(i\in I\) greater than all 
    \(j\in J\). Since our construction will ensure that 
    \(|J|\leq\kappa\) and the partial order \(I\) (or domain category)
    is \(\kappa\)-directed (-filtered) -- we will have at least one 
    such \(i\) and we will have factored \(f\) appropriately.
\end{lemma}
\begin{proof}
    The reader must forgive the following proof, as it isn't quite as 
    insightful as to move or impart the reader with any deep beauty,
    but it is baroque enough to possibly deeply confuse them. 

    Take some \(f:X\to\colim_kD_k\) for the rest of the proof.
    
    \paragraph{Defining \(J\):}
    Given \(x,y\in X\), consider the set 
    \[
        \Delta(x,y) = \set{\delta_i(a,b)}{
            i\in I
            \quad
            a, b\in D_i
            \quad
            \coproj_i(a) = f(x)
            \quad
            \coproj_i(b) = f(y)
        }
    \]
    (Note that \(\Delta(x,y)\) is no emptier than \(X\), as 
    \(f(x)\) must “be” (as in, up to \(\coproj\)) in some 
    \(D_i\), and \(f(y)\), in some \(D_j\). And
    thus there is some \(k\) greater than both must “be” in at
    the same time.)
    
    If we admit the axiom of choice (which we do, and must to if 
    we even want to start talking about things like 
    \(|\quantale Q|\) for non-special quantales), we can choose a 
    a \(\Xi(x,y)\)
    \[
        \Xi(x,y) = \set{(i, a, b)}{
            i\in I
            \quad
            a, b\in D_i
            \quad
            \coproj_i(a) = f(x)
            \quad
            \coproj_i(b) = f(y)
        }
    \]
    such that \(\Xi(x,y)\iso\Delta(x,y)\) and the isomorphism 
    is given by the map 
    \[
        (i,a,b)\mapsto\delta_i(a,b)
    \]
    And now, if we take the projection of \(\Xi(x,y)\) into the 
    first coordinate, corresponding to the index \(i\) that
    \(a,b\) inhabit, we get a set we call \(\Gamma(x,y)\):
    \[
        \Gamma(x,y) = \pi_I[\Xi(x,y)]
    \]

    Now, the reader may be assured: if anything \(\Gamma(x,y)\) 
    must have cardinalty stricly less than \(\kappa\). As 
    \(\Delta(x,y)\) is a subset of \(\quantale Q\), which is 
    itself smaller than \(\kappa\) -- and we have only ever 
    decreased its size by applying the above constructions.
    And so, we may indeed proceed with defining \(J\) 
    \[
        J = \bigcup_{x,y}\Gamma(x,y)
    \]
    Since \(|X|<\kappa\) so does \(|X\times X|\) (provided, 
    say, \(\kappa\) is infinite. We'd force \(\kappa\geq\omega\)
    otherwise, no harm done. 

    Since \(|X\times X|<\kappa\) it follows that we are doing
    a union of less than \(\kappa\) sets of cardinalities that 
    are lesser than \(\kappa\) -- and regularity \emph{is} the
    fact that that itself must be smaller than \(\kappa\).
    And hence \(|J|\) is in fact still smaller than \(\kappa\).
    
    Let, therefore, \(\gamma\in I\) be some element greater 
    than all \(j\in J\). Obviously the same can be done for 
    filtered diagrams as opposed to posets. Again: no harm is
    done to the argument.
    
    \paragraph{Factoring \(f\) through \(D_\gamma\):}
    Recall the construction for colimits, we may regard 
    colimits as appropriate quotients of disjoints unions, 
    this is now be centrally useful:
    
    Given \(x\in X\), its image under \(f\) is an equivalence
    class \(\eqclass{(a, i)}\) for an \(a\in D_i\) an elected
    representative (not necessarily democratically, the Axiom 
    of Choice doesn't imply that everyone necessarily has a 
    say). Hence, it is evident that by taking the indices of 
    those representatives we shall have a set \(J'\) of 
    cardinalty smaller than \(\kappa\), to which is associated
    an element \(\gamma'\)---greater than all its elements---%
    such that we have a function 
    \(X\xrightarrow{\bar f} D_{\gamma'}\) taking \(x\) to 
    \(D(i\leq\gamma')(a)\), which factors \(f\).

    \paragraph{Betrayal!}
    However, this function is \emph{not} necessarily a 
    functional \emph{morphism} because we aren't taking 
    \(\delta\) into account! \(\delta(x,y)\) might not be less 
    than \(\delta(\bar f(x),\bar f(y))\).
    
    It \emph{obviously} factors \(f\)... in \(\Set\). 
    Despite much fear and trembling (some loathing too), or 
    largely because of them, it is possible to enhance it into 
    a \emph{morphism} that factors it. This is done with the 
    help of our old friend \(\gamma\).
    
    \paragraph{Redemption:}
    A convenient fact left out (for dramatic purposes) of the
    first part of the proof is that 
    \[
        \delta(f(x),f(y)) = \sup \Delta(x,y)
    \]
    Which is because how \(\delta\) is defined on colimits:
    \[
        \delta(\eqclass a,\eqclass b) = 
        \bigvee_{i\in I}\bigvee_{\substack{
            \alpha\in\eqclass a\cap D_i\\
            \beta \in\eqclass b\cap D_i
        }}\delta_i(\alpha,\beta)
    \]
    But also recall that \(D_\gamma\) is “above” all \(D_j\)
    for \(j\in J\):
    \[
        \delta_j(u,v)\leq
        \delta_\gamma(
            D(j\leq\gamma)(u),
            D(j\leq\gamma)(v)
        )
    \]
    Hence, if \(\coproj_j(a) = f(x)\) and \(\coproj_j(b) = f(y)\)
    for some \(j\), it must follow that
    \[
        \delta_j(a,b)\leq
        \delta(
            D(j\leq\gamma)(a),
            D(j\leq\gamma)(b)
        )
    \]
    But taking the supremum of the \((i,a,b)\) that do that is 
    just \(\delta(f(x),g(y)\). We also know that it cannot 
    grow any more than that and so we obtain:
    \[
        \delta(
            D(j\leq\gamma)(a),
            D(j\leq\gamma)(b)
        ) = 
        \delta(f(x),f(y))
    \]

    This tells us that \(\gamma\) can do all that \(\gamma'\) 
    could -- in that we can find preimages of \(f(x)\) under 
    \(\coproj\) for every \(x\) -- but also we can do so
    in such a way that is a \emph{morphism}. We only need to 
    concerns ourselves with \(\delta\), as the extent is 
    trivially always preserved. This means we have indeed 
    factored \(f\) as 
    \[
        X
        \xrightarrow{\bar f} 
        D_\gamma
        \xrightarrow{\coproj_\gamma}
        \colim\limits_kD_k
    \]
    and that tells us:
    \[
        \phi(\eqclass{\bar f}) = 
        \coproj_\gamma\circ\bar f = f
    \]
\end{proof}

\begin{lemma}[\(\phi\) is injective]
    The converse of the above holds as well:
    \(\phi\) is injective.
\end{lemma}
\begin{proof}
    Take \(X\) and \(D\) as above, we ought to show that 
    \(\phi\) is injective and thus that if 
    \(\phi(\eqclass{(f, i)}) = \phi(\eqclass{(g, j)})\)
    we ought to be able to show that 
    \(\eqclass{(f, i)} = \eqclass{(g, j)}\). 
    That is bound to be fun. Suppose, then, we do have 
    such a pair. It follows definitionally that 
    \[
        \coproj_i\circ f = \coproj_j\circ g
    \]
    Since these functions are extensionally the same, we have 
    that for each \(x\in X\), \(f(x)\sim g(x)\). Where this 
    equivalence is the equivalence the symmetric transitive 
    closure of
    \[
        (a\in D_i)\sim([D(i\leq j)](a)\in D_j)
    \]
    defined on \(\coprod_i D_i\). Hence, it amounts to saying
    that \((a,i)\sim(b,j)\iff\) there is a messy zig-zag 
    diagram, as below, connecting them
    \\[1em]
    \adjustbox{scale=0.7,center}
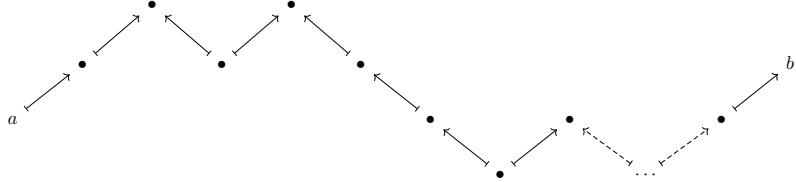
{\begin{tikzcd}
    	&& \bullet && \bullet \\
    	& \bullet && \bullet && \bullet &&&&&& b\\
    	a &&&&&& \bullet && \bullet && \bullet \\
    	&&&&&&& \bullet && \cdots
    	\arrow[maps to, from=3-1, to=2-2]
    	\arrow[maps to, from=2-2, to=1-3]
    	\arrow[maps to, from=2-4, to=1-3]
    	\arrow[maps to, from=2-4, to=1-5]
    	\arrow[maps to, from=2-6, to=1-5]
    	\arrow[maps to, from=3-7, to=2-6]
    	\arrow[maps to, from=4-8, to=3-7]
    	\arrow[maps to, from=4-8, to=3-9]
    	\arrow[maps to, from=3-11, to=2-12]
    	\arrow[dashed, maps to, from=4-10, to=3-9]
    	\arrow[dashed, maps to, from=4-10, to=3-11]
    \end{tikzcd}}
    \\[1em]
    Consequently, to each \(x\) there is (at least one)
    finite diagram as above connecting \(f(x)\) and 
    \(g(x)\). Each of those elected diagrams concerns 
    only finitely many \(i\) in \(I\) -- and there are 
    less than \(\kappa\) \(x\) in \(X\). And thus if we 
    take them all together, we will still have what 
    amounts to less than \(\kappa\) indices. We call one 
    such collection of indices simply \(J\), nevermind 
    which one exactly.

    We take some \(\gamma\) greater than all \(j\) in 
    \(J\), and once again consider \(D_\gamma\). Since
    The diagram formed by 
    \(D\restriction(J\mathbin{\dot\cup}\{\gamma\})\) is
    commutative, we have embedded the messy zig-zag
    for \emph{every} \(x\) in a beautifully commutative 
    way:
    
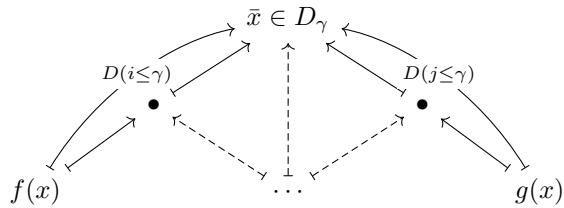
\[\begin{tikzcd}
    	&& {\bar x\in D_\gamma} \\
    	& \bullet && \bullet \\
    	{f(x)} && \cdots && {g(x)}
    	\arrow[dashed, maps to, from=3-3, to=2-4]
    	\arrow[maps to, from=3-5, to=2-4]
    	\arrow[dashed, maps to, from=3-3, to=2-2]
    	\arrow[maps to, from=3-1, to=2-2]
    	\arrow["{D(i\leq \gamma)}"{description}, curve={height=-18pt}, maps to, from=3-1, to=1-3]
    	\arrow[maps to, from=2-2, to=1-3]
    	\arrow[dashed, maps to, from=3-3, to=1-3]
    	\arrow[maps to, from=2-4, to=1-3]
    	\arrow["{D(j\leq \gamma)}"{description}, curve={height=18pt}, maps to, from=3-5, to=1-3]
    \end{tikzcd}\]
    Hence, \(g(x)\) and \(f(x)\) get identified 
    \emph{for every \(x\) at the same time} in this
    particular \(D_\gamma\). In particular, we may 
    simply take \(\bar f\) to be the morphism taking 
    \(x\) to \(\bar x\), which coincides with 
    \(\bar g\) -- doing the same. But
    \begin{align*}
        \hom(X,D(j\leq \gamma))(g) 
        &=
        [D(g\leq \gamma)] \circ g
        \\&=
        \bar g 
        \\&=
        \bar f 
        \\&=
        [D(i\leq \gamma)] \circ f 
        \\&=
        \hom(X,D(i\leq \gamma))(f) 
    \end{align*}
    And hence, we have a messy zig-zag (this time in \(\Set\))
    connecting those arrows, and hence in the colimit they 
    are actually the same:
    \[
        \eqclass{(f,i)} = \eqclass{(g,j)}
    \]
    and this is what we set out to prove.
\end{proof}

\begin{theorem}[\(|X|<\kappa\) implies \(\kappa\)-compactness]
    Which is obvious in the light of the lemmas above.
\end{theorem}

\subsection{Accessibility and Presentability}

\begin{theorem}[\(\QSet\) is indeed \(\kappa\)-accessible]
    “Trivial”.
\end{theorem}
\begin{proof}
    Unsurprisingly, if you take \(Y\) some \(\quantale Q\)-set
    and take the inclusion poset of \(\parts^{<\kappa}(|Y|)\)
    which is always (due to regularity) \(\kappa\)-directed.

    In more precise terms, those “small” parts of \(Y\) get 
    mapped to... themselves, with \(\delta\) given by 
    restriction. It is immediately evident that the colimit 
    of this diagram is \(Y\): 

    Take some \(y\in Y\), it corresponds to the equivalence 
    class containing manifold copies of itself but tagged 
    with whichever subset that happened to contribute its 
    membership to the disjoint union.

    Take a pair \(y,y'\), since \(\delta\) on subsets will
    be given by restriction, the \(\delta\) of their 
    corresponding classes will just be their \(\delta\). 
    So the colimit is evidently just isomorphic to \(Y\).
\end{proof}

\begin{theorem}[Local Presentability]
    Since we already know that \(\QSet\) is a cocomplete 
    category, we have actually shown that it is locally
    presentable. Since accessible cocomplete categories are
    invariably locally presentable.
\end{theorem}
	\section{Monoidal Structures}

It is already known, thanks to the construction of limits in a previous
section, that \(\QSet\) has a “reasonable” monoidal category structure.
This product, however, doesn't often have an exponential associated to 
it: if it did, \(\blank\times X\) would be cocomplete. Recall the 
construction of the categorical products:
\[
    X\times Y = \set{
        (x,y)
    }{
        \extent_X x = \extent_Y y
    }
\]
with \(\delta((x,y),(a,b)) = \delta_X(x,a)\wedge\delta(y,b)\).
Since coequalizers have to do with taking suprema, the product being
cocomplete would mean that something like
\[
    a \wedge \bigvee_i b_i = \bigvee_i a\wedge b_i
\]
would have to hold. This would make \(\quantale Q\)'s underlying 
lattice a locale---although it does not force \(\tensor=\wedge\), of
course.

The question the arises: are there monoidal closed (and semicartesian)
category structures naturally defined over \(\QSet\) and how are those 
different structures related to each other? We offer some such 
structures -- arranged in a hierarchy of monoidal products related to 
\(\tensor\).

From the product construction, we can spot two extension points: 
we can change \(\delta\) to use \(\tensor\) as opposed to \(\wedge\)
and we can use some other relation \(\sim\) between \(\extent_X\) and
\(\extent_Y\) as opposed to \(=\). The former we always take, the 
latter requires some consideration so that desirable categorical 
properties may still hold.

\begin{definition}[Locallic Congruence]
    A locallic congruence is an equivalence relation on a locale 
    such that 
    \[
        [\forall i\in I : a\sim b_i]\implies a\sim\bigvee_{i\in I}b_i
    \]
    \[
        [a\sim b, a\sim c] \implies [a\sim (b\wedge c)]
    \]
    We say “a locallic congruence over a \(\quantale Q\)” for a 
    quantale \(\quantale Q\) meaning one such congruence over its 
    locale of idempotent elements.
\end{definition}

\begin{definition}[Congruential Tensor]
    Our tensors come from locallic congruences on \(\quantale Q\).
    Namely, taking one such congruence \(\sim\), we define the 
    operation
    \[
        X\tensor Y = \set{(x,y)}{\extent_Xx\sim\extent_Yy}
    \]
    \[
        \delta((x,y),(a,b)) = \delta_X(x,a)\tensor\delta(y,b)
    \]
    We should also define its action on morphisms: which is to take 
    \((f:X\to Y,~g:A\to B)\) to 
    \[
        (x,a) \xmapsto{f\tensor g} (f(x),g(a))
    \]
    Functoriality is trivial.
    
    We claim the above defines an obvious functor, that this is 
    actually associative and commutative, that it is semicartesian,
    cocomplete, and that it has a unit. We shall first provide the 
    definition for the unit, and then proceed with the appropriate 
    proofs for our claims.
\end{definition}

\begin{definition}[Congruential Tensor Unit]
    Given the above \(\sim\), it defines equivalence classes on 
    \(\extent\quantale Q\), and since \(\sim\) is closed under 
    suprema:
    \[
        a\sim\sup\eqclass a
    \]
    And so, the set \(\extent\quantale Q\) is in bijection with the
    following regular subterminal, given by 
    \[
        \set{\sup\eqclass a}{a\in\extent\quantale Q}
    \]
    which is what we take to be \(1\) -- the claimed unit for 
    \(\tensor\).
\end{definition}

\begin{remark}
    In this section, we shall -- for notational reasons -- be using
    \(\tensor\) as a functor to denote a generic congruential 
    tensor -- coming from some fixed but generic \(\sim\). Where
    necessary/convenient, we may specify the relation to 
    disambiguate. Later own, though, we shall use \(\tensor\) 
    referring to the minimal tensor, given by \(a\sim b\iff a=b\).
    Similarly, we shall refer to the maximal tensor, given by the 
    chaotic relation by the symbol \(\chaosor\).
\end{remark}

\begin{lemma}[\(\tensor\) is commutative]
\end{lemma}
\begin{proof}
    Since \(\tensor:\quantale Q\times\quantale Q\to\quantale Q\)%
    ---the actual algebraic operation on the quantale---is taken 
    to be commutative, and the fibration over the extents is over
    a symmetric relation, it is obvious that it will be 
    commutative.
\end{proof}

\begin{lemma}[\(\tensor\) is associative]
\end{lemma}
\begin{proof}
    Associativity has a canonical isomorphism we ought to consider
    -- that being what we would do in \(\Set\) if we were given the
    task there: \((a,(b,c))\mapsto((a,b),c)\). We ought to show 
    that this is an isomorphism, instead of just a function between
    the underlying sets. This is achieved by realizing that
    \[
        (a,(b,c))\in A\tensor(B\tensor C) 
        \iff 
        \extent a\sim\extent b\sim\extent c
        \iff
        ((a,b),c)\in (A\tensor B)\tensor C
    \]
    And we do \emph{that} using \(\wedge\)-congruence: 
    \begin{align*}
            (a,(b,c))\in A\tensor(B\tensor C) 
        &\implies
            [\extent b \sim \extent c] 
            ~\text{and}~
            [\extent a \sim \extent b\tensor\extent c]
        \\&\implies
            [\extent b \sim \extent b\tensor\extent c] 
            ~\text{and}~
            [\extent a \sim \extent b\tensor\extent c]
        \\&\implies
            [\extent a\sim\extent b\sim\extent c]
    \end{align*}
    Obviously, mutatis mutandis, one can prove the same for 
    \(((a,b),c)\in(A\tensor B)\tensor C\) -- this can also 
    be seen as a consequence of the above and commutativity.
    Obviously, if the extent equivalence chain holds, then
    one can form the triples in either shape, so the 
    logical equivalence holds.

    We have shown that the obvious associator is in fact a 
    function, but it is also bound to be an isomorphism 
    since it is evidently a bijection and preserves 
    \(\delta\) thanks to the associativity of \(\tensor\)
    as a \(\quantale Q\) operation.

    We call the above associator “\(\alpha\)”. It is easy
    to see that \(\alpha\) is natural and indeed satisfies
    the pentagon identity. For the skeptical readers, a 
    proof can be sketched: the forgetful functor back into
    \(\Set\) does not make commutative any diagram that 
    didn't already enjoy the property -- and the forgetful 
    image of our pentagon identity is the restriction of a 
    commutative diagram in set (the pentagon identity for 
    the set-theoretical cartesian product).
\end{proof}

\begin{lemma}[\(\tensor\) has \(1\) as a unit]
\end{lemma}
\begin{proof}
    We haven't given \(1\)'s \(\delta\) but have left it 
    implicit in saying it is a regular subterminal -- meaning
    it is a subset of \(\extent Q\) with \(\delta\) defined as
    \(\wedge\). This suffices to show that it is indeed a 
    \(\quantale Q\)-set.

    Suppose now we have some \(\quantale Q\)-set \(X\) and 
    let us consider its product with \(1\): 
    \[
        X\tensor 1 = \set{(x,\sup\eqclass e)}{
            e\in\extent\quantale Q,~
            \extent_Xx \sim\extent_1\sup\eqclass e
        }
    \]
    And hence, we know that \(\extent_Xx\sim e\) if and
    only if \((x,\sup\eqclass e)\in X\tensor 1\). Since 
    \(\sim\) is symmetric, it follows that the only element
    that can ever get paired with \(x\) is 
    \(\sup\eqclass{\extent x}\). And hence \(|X|\) is in 
    bijection with \(|X\tensor 1|\). This bijection is 
    naturally \(\delta\)-preserving:
    \begin{align*}
        \delta(x,y) 
        &=  \delta(x,y)\tensor\extent x\tensor\extent y
        \\&=\delta(x,y)\tensor\extent x\tensor\extent y
                \tensor\sup\eqclass{\extent x}
                \tensor\sup\eqclass{\extent y}
        \\&=\delta(x,y)
                \tensor\sup\eqclass{\extent x}
                \tensor\sup\eqclass{\extent y}
        \\&=\delta(x,y)\tensor\delta_1(
            \sup\eqclass{\extent x},
            \sup\eqclass{\extent y}
        )
        \\&=\delta(
            (x,\sup\eqclass{\extent x}), 
            (y,\sup\eqclass{\extent y})
        )
    \end{align*}

    One has projections \(X\tensor 1, 1\tensor X\to X\)
    given by forgetting the second coordinate (or just 
    taking the inverse of the bijection we have 
    established above). Those are what will become our 
    unitors: The morphisms \(\rho:X\tensor 1\to X\) and 
    \(\lambda:1\tensor X\to X\) are both obviously natural 
    in \(X\) and we won't spend any time proving it here. 
\end{proof}

\begin{theorem}[\((\QSet, \tensor, 1)\) is a symmetric monoidal category]
    In light of what we've established, we would have to show
    that the associator and unitor satisfy the triangle identity 
    and that we have a braiding satisfying one hexagonal identity
    and the symmetry condition.
\end{theorem}
\begin{proof}
    We establish the unitor-associator triangle identity by 
    tracing around an element around the path, noting that 
    \(\extent(x,\sup\eqclass{\extent x}) = \extent x\) and 
    hence \(\extent x\sim\extent y\).
    \[\begin{tikzcd}
    	{(X\tensor 1)\tensor Y} && {X\tensor (1\tensor Y)} \\
    	& {X\tensor Y} \\
    	{((x,\sup\eqclass{\extent x}),y)} && {(x,(\sup\eqclass{\extent x},y))} \\
    	\\
    	{(x,y)} && {(x,(\sup\eqclass{\extent y},y))}
    	\arrow["\alpha", from=1-1, to=1-3]
    	\arrow["{\rho\tensor Y}"{description}, from=1-1, to=2-2]
    	\arrow["X\tensor\lambda"{description}, from=1-3, to=2-2]
    	\arrow["{\text{\(\extent x\sim\extent y\)}}", Rightarrow, no head, from=3-3, to=5-3]
    	\arrow["X\tensor\lambda", maps to, from=5-3, to=5-1]
    	\arrow["\alpha", maps to, from=3-1, to=3-3]
    	\arrow["{\rho\tensor Y}"', maps to, from=3-1, to=5-1]
    \end{tikzcd}\]

    The braiding in question is simple: \((x,y)\mapsto(y,x)\) 
    which can be formed as \(\sim\) is symmetric. Moreover, 
    it is obviously \(\delta\)-preserving and bijective, making 
    it an isomorphism. It is trivially natural on \(X\) and \(Y\), 
    and swapping \(X\) for \(Y\) and vice versa evidently shows 
    that the symmetry condition holds. So all that remains is to 
    show the hexagon braiding identity. For now, let's give a name
    to the braiding: \(\beta : X\tensor Y\to Y\tensor Y\):
    \[\begin{tikzcd}
    	& {X\tensor (Y\tensor Z)} & {(Y\tensor Z)\tensor X} \\
    	{(X\tensor Y)\tensor Z} &&& {Y\tensor (Z\tensor X)} \\
    	& {(Y\tensor X)\tensor Z} & {Y\tensor(X\tensor Z)}
    	\arrow["\alpha", from=2-1, to=1-2]
    	\arrow["\beta", from=1-2, to=1-3]
    	\arrow["\alpha", from=1-3, to=2-4]
    	\arrow["{\beta\tensor Z}"', from=2-1, to=3-2]
    	\arrow["\alpha"', from=3-2, to=3-3]
    	\arrow["Y\tensor\beta"', from=3-3, to=2-4]
    \end{tikzcd}\]
    Here, the reader is invited to trace an element's orbit along
    those paths, there are no traps -- we swear.
\end{proof}

\begin{theorem}[\(\tensor\) is cocomplete in either entry]
\end{theorem}
\begin{proof}
    Showing it directly may be painful, so instead we opt to do 
    it in two steps: proving it preserves coproducts and 
    proving it preserves coequalizers. Since all colimits are 
    coequalizers of coproducts, preserving both means to preserve
    all.

    \paragraph{Coproduct Preservation}
    Take \(X_i\) for \(i\in I\) an indexed set of \(\quantale Q\)%
    -sets. We can inspect \(\left(\coprod_i X_i\right)\tensor Y\)
    and its elements are 
    \[
        ((x,i),y)~\text{\(x\in X_i\) st. \(\extent_i x \sim y\)}
    \]
    Inspecting the elements of \(\coprod_i(X_i\tensor Y)\) yields
    \[
        ((x,y),i)~\text{\(x\in X_i\) st. \(\extent_i x \sim y\)}
    \]
    So it's just a shuffling of the label saying which index that
    instance of \(x\) belongs to. Completely immaterial. Their 
    \(\delta\)s just as similar:
    \begin{align*}
        \delta(((x,i),y), ((x',i'),y')) &= 
        \delta((x,i),(x',i'))\tensor\delta(y,y') 
        \\&= 
        \begin{cases}
            \delta(y,y')\tensor\delta_i(x, x'), & i = i',\\
            \delta(y,y')\tensor\bot,            & i \not= i'.
        \end{cases}
        \\
        \\
        \delta(((x,y),i), ((x',y'),i')) &= 
        \begin{cases}
            \delta_{X_i\tensor Y}((x,y),(x',y')), & i = i',\\
            \bot,                                 & i \not= i'.
        \end{cases} \\&= 
        \begin{cases}
            \delta_i(x,x')\tensor\delta(y,y')  & i = i',\\
            \bot,                              & i \not= i'.
        \end{cases}
    \end{align*}
    Since \(\bot\) is absorbing, they are all the same.

    \paragraph{Coequalizer Preservation}
    Take \(f_i:A\to B\) for a pair of fixed \(A\) and \(B\),
    let's consider the coequalizer
    \[\begin{tikzcd}
    	A && B & C
    	\arrow[""{name=0, anchor=center, inner sep=0}, "{f_i}"{description}, curve={height=-18pt}, from=1-1, to=1-3]
    	\arrow[""{name=1, anchor=center, inner sep=0}, "{f_j}"{description}, curve={height=18pt}, from=1-1, to=1-3]
    	\arrow["{\eqclass{\blank}}", from=1-3, to=1-4]
    	\arrow["\vdots"{description}, draw=none, from=0, to=1]
    \end{tikzcd}\]
    which, as we know, is given by the quotient of the 
    reflexive and transitive closure of the relation given 
    by
    \[
        f_i(a) \approx f_j(a)
    \]

    So consider \(X\tensor f_i:X\tensor A\to X\tensor B\). 
    First we note that an element of their coequalizer 
    must be some \(\eqclass{(x,b)}_{\approx'}\) with 
    \(\approx'\) generated by 
    \[
        (X\tensor f_i)(x,a) = (x,f_i(a)) \approx'
        (x,f_j(a)) = (X\tensor f_j)(x,a)
    \]
    hence it is obvious (since it's the same \(x\) throughout 
    the line and \(\approx'\) is generated freely by the above
    conditions, that an equivalence class of the relation as 
    defined above will and can only ever have the same \(x\) 
    in all of its members. 

    consider the projection on the second coordinate, ie. 
    \(\pi_B\eqclass{(x,b)}\) as a subset of \(B\). It is 
    evident that all \(b'\) in such an assembly will be such
    that \(b'\approx b\) -- by definition. It is also true 
    that \(\extent b' = \extent b\) since the equivalence 
    relation requires them to have come from the same point
    by possibly different -- and zigzagging (but that doesn't
    matter, since extent will be preserved) -- paths. 

    Since the same will hold for the broader subset that is
    the equivalence class \(\eqclass{b}_\approx\), we can 
    confidently say that 
    \(\extent x\sim\extent\eqclass{b}_\approx\). Thus, we 
    can safely hoist the \(X\) component out of 
    each equivalence class \(\eqclass{(x,b)}\) and obtain an
    element \((x,\eqclass{b}_\approx)\in X\tensor C\).

    Similarly, any such element will necessarily be such 
    that we can form \((x,b')\) for every 
    \(b'\in\eqclass{b}_\approx\) and hence comes from an 
    element \(\eqclass{(x,b)}_{\approx'}\) in the manner 
    described above. This surjection we have shown is 
    obviously also injective. So for an isomorphism, all
    that we would have to show is \(\delta\)-preservation.
    \begin{align*}
        \delta((x,\eqclass b), (y,\eqclass\beta))
        &=  \delta(x,y)\tensor\delta(\eqclass b,\eqclass\beta)
        \\&=\delta(x,y)\tensor\bigvee_{\substack{
            u\in\eqclass b\\
            v\in\eqclass\beta
        }}
            \delta(u,v)
        \\&=\bigvee_{\substack{
            u\in\eqclass b\\
            v\in\eqclass\beta
        }}
            \delta(x,y)\tensor\delta(u,v)
        \\&=\bigvee_{\substack{
            u\in\eqclass b\\
            v\in\eqclass\beta
        }}
            \delta((x,u),(y,v))
        \\&=\delta(\eqclass{(x,u)},\eqclass{(y,v)})
    \end{align*}

\end{proof}

\begin{theorem}[\((\QSet,\tensor)\) is monoidal closed]
\end{theorem}
\begin{proof}
        \(X\tensor\blank\) is a cocomplete
        endofunctor on a locally presentable category.
\end{proof}

\subsection{Formalizing their Hierarchy}
As previously seen, an equivalence relation \(\sim\) 
gives rise to a tensor \(\tensor[\sim]\) which we had 
heretofore neglected to mark with the relation that 
originated them. Since relations are ordered by 
inclusion (and hence, by implication) -- it would be 
strange if some kind of similar hierarchy did not 
connect their tensorial spawn.

To that effect, we have introduced a notion of 
morphism between monoidal categories called “translax”
monoidal functors -- which are neither lax nor oplax 
functors but do indeed form a 2-category with objects
being monoidal categories.

\begin{definition}[Translax Monoidal Functor]
    Roughly speaking, a translax monoidal functor is 
    a akin to a monoidal functor but the arrows 
    associated to units are going in the reverse 
    (hence trans) direction to the monoidal product 
    side of things.

    In more formal terms, given monoidal categories
    \(\category A\) and \(\category B\), a monoidal
    functor \(F:\category A\to\category B\) is 
    a functor endowed with a certain tensor 
    “covariant” transformation and a “contravariant”
    unit map satisfying some coherence conditions. 
    This is to say, a functor \(F\) between the 
    underlying categories and
    \begin{align*}
        \mu: F(\blank \tensor_\category A  \blank)
        &\to F(\blank)\tensor_\category BF(\blank)\\
        \epsilon: 1_\category B&\to F(1_\category A)
    \end{align*}
    satisfying\\
    \adjustbox{scale=0.75,center}{\begin{tikzcd}
    	&& {F((X\tensor Y)\tensor Z)} \\
    	{F(X\tensor(Y\tensor Z))} &&&& {F(X\tensor Y)\tensor F(Z)} \\
    	\\
    	{F(X)\tensor F(Y\tensor Z)} &&&& {(F(X)\tensor F(Y))\tensor F(Z)} \\
    	&& {F(X)\tensor (F(Y)\tensor F(Z))}
    	\arrow["{F(\alpha)}"{description}, from=1-3, to=2-1]
    	\arrow["\mu"{description}, from=1-3, to=2-5]
    	\arrow["{\mu\tensor F(Z)}"{description}, from=2-5, to=4-5]
    	\arrow["\mu"{description}, from=2-1, to=4-1]
    	\arrow["{F(X)\tensor\mu}"{description}, from=4-1, to=5-3]
    	\arrow["\alpha"{description}, from=4-5, to=5-3]
    \end{tikzcd}}
    \[\begin{tikzcd}
    	{F(X)} & {1\tensor F(X)} && {F(X)\tensor 1} & {F(X)} \\
    	{F(1\tensor F(X))} & {F(1)\tensor F(X)} && {F(X)\tensor F(1)} & {F(X\tensor 1)}
    	\arrow["{F(\rho^{-1})}"', from=1-5, to=2-5]
    	\arrow[from=2-5, to=2-4]
    	\arrow["{\rho^{-1}}"{description}, from=1-5, to=1-4]
    	\arrow["{F(X)\tensor\epsilon}"', from=1-4, to=2-4]
    	\arrow["{\epsilon\tensor F(X)}", from=1-2, to=2-2]
    	\arrow[from=2-1, to=2-2]
    	\arrow["{F(\lambda^{-1})}"', from=1-1, to=2-1]
    	\arrow["{\lambda^{-1}}", from=1-1, to=1-2]
    \end{tikzcd}\]
\end{definition}
\begin{theorem}
    The obvious composition makes this a category, with the 
    natural identity.
\end{theorem}
\begin{proof}
\end{proof}

Here we now claim that, in fact, our construction
\({\sim}\mapsto\tensor[\sim]\)
is actually a functor from the lattice of locallic 
congruences over \(\extent\quantale Q\) to the 2-category of
monoidal categories and translax functors between them.
In our case, all components of the involved translax 
functors are all \emph{inclusions} -- regular monos.

\begin{lemma}[Unit Inclusion]
    Suppose \(\sim\implies\approx\), there is a regular mono
    \(\epsilon:1_\approx\to1_\sim\) that is 
    given functorially with respect to \(\leq\).
\end{lemma}
\begin{proof}
    Take an element \(b = \sup\eqclass{b}_\approx\), since 
    \(a\sim b\implies a\approx b\), we have that \(b\) is 
    also \(\sup\eqclass{b}_\sim\). This simple observation,
    that laxer equivalences have fewer equivalent classes 
    partitioning totality, also means that our elected 
    representatives of those classes will exist in 
    equivalences that \emph{refine} ours.

    Hence, there is an inclusion from the unit of the 
    larger tensor to the unit of the smaller one. This is 
    obviously (contravariantly) functorial, as they are 
    just inclusions.
\end{proof}

\begin{lemma}[Product Inclusion]
    There is a functorial assignment of natural 
    transformations over the set of locallic congruences 
    over \(\quantale Q\).
\end{lemma}
\begin{proof}
    Since monoidal products are defined through set 
    comprehension, it is evident that 
    \[
        [\sim\implies\approx]\implies[
            (x,y)\in X\tensor[\sim]    Y \implies
            (x,y)\in X\tensor[\approx] Y
        ]
    \]
    Consequently, we have maps -- that are trivially 
    natural on \(X\) and \(Y\) -- given “functorially”.
    This is to say: if we take the category of functors
    \(\QSet\times\QSet\to\QSet\) and natural 
    transformations between them, The mapping taking 
    \(\sim\) to \(\tensor[\sim]\) and taking 
    \({\sim}\leq{\approx}\) to  
    \(\tensor[\sim]\hookrightarrow\tensor[\approx]\)
    is trivially functorial, as the maps are just 
    inclusions.
\end{proof}

\begin{theorem}
    There is a functor from the implication category of 
    locallic congrunces over \(\quantale Q\) and the 
    2-category of monoidal categories with translax 
    functors. That functor takes \(\sim\) to its 
    associated tensor product.
\end{theorem}
\begin{proof}
    We already know that the product inclusion is a 
    natural transformation -- functorially dependent
    on \(\sim\) -- of the appropriate type. And we know
    that the unit inclusions are contravariantly 
    functorial with respect to \(\sim\). Therefore, the
    obvious choice of \(F\) is \(\id_{\QSet}\). What 
    remains is to show that the above choices jointly 
    form a translax functor.

    This verification, however, is rather dull to read
    and to transcribe from our notes -- we shall omit it
    since it is quite trivial.
\end{proof}
	\section{Change of basis}

There are plenty of possible definitions for morphisms between
quantales. The most basic that is remotely useful is to be an
order morphism jointly with a semigroup morphism. There are many
examples of such morphisms, such as subquantale inclusions, 
projections.

There are nontrivial such morphisms, such as when \(\quantale Q\)
is semicartesian and commutative, in which one may form 
\[
    q \mapsto q^- = \max\set{e\in\extent Q}{e\preceq q}
\]
which happens to be right adjoint to the inclusion of idempotents
into the quantale. This means that relations arising from such 
morphisms would relate locallic-sets (which are topos-adjacent)
and our \(\quantale Q\)-sets. Thus, for now, let us explore these 
simple morphisms which are just non-decreasing functions that 
preserve products, and the obvious functor that arises from them.

\begin{definition}
    Given \(f:\quantale P\to\quantale Q\), one defines 
    \(f_*:\QSet[\quantale P]\to\QSet[\quantale Q]\) to be the 
    functor
    \[
        f_*(X,\delta) = (X,f\circ\delta)
    \]
    with the trivial action on morphisms. As defined, this 
    obviously preserves the identity. Moreover, its basically 
    immediate that -- if \((X,f\circ\delta)\) is indeed always 
    a \(\quantale Q\)-set whenever \((X,\delta)\) is a 
    \(\quantale P\)-set -- that the action if functorial: we are
    just composing functions after all.

    To see the claim that remains, simply realize that all 
    \(\quantale Q\)-set axioms are either trivially valid 
    (symmetry) or depend on \(\leq\) and \(\tensor\) to be 
    preserved.
\end{definition}

    \printbibliography{}
 
\end{document}